\documentclass[11pt]{amsart}
\usepackage{amssymb}
\usepackage{amsmath}
\usepackage{amsthm}

\newtheorem{thm}{Theorem}[section]
\newtheorem{lem}[thm]{Lemma}

\newtheorem{prop}[thm]{Proposition}
\theoremstyle{definition}

\newtheorem{rem}[thm]{Remark}

\newcommand{\R}{\mathbb{R}}

\begin{document}
\title[Finite Morse index solutions and asymptotics]{Finite Morse index solutions and asymptotics of weighted nonlinear elliptic equations$^*$}
\thanks{}
\author[Y. Du and Z.M. Guo]{Yihong Du$^\dag$ and Zongming Guo$^\ddag$}
\thanks{$^\dag$ School of Science and Technology,
University of New England, Armidale, NSW 2351, Australia
(ydu@\allowbreak turing.\allowbreak une.\allowbreak edu.\allowbreak
au).}
\thanks{$^\ddag$ Department of Mathematics, Henan Normal University,
Xinxiang, 453007, P.R. China (gzm@htu.cn)}

\subjclass{Primary 35B45; Secondary 35J40}

\thanks{* The research of Y. Du was supported by the Australian Research Council, and that of Z.M. Guo was supported by NSFC
(11171092, 10871060) and Innovation Scientists and Technicians Troop
Construction Projects of Henan Province (114200510011).}

\keywords{Nonlinear Schr\"odinger equations with Hardy potential,
supercritical exponent, stability, finite Morse index solutions,
asymptotics}
\date{\today}
\def\baselinestretch{1}

\begin{abstract}
By introducing a suitable setting, we study the behavior of finite
Morse index solutions of the equation
\[
-\mbox{div} (|x|^\theta \nabla v)=|x|^l |v|^{p-1}v \;\;\; \mbox{in
$\Omega \subset \R^N \; (N \geq 2)$}, \leqno(1)
\]
where $p>1$, $\theta, l\in\R^1$ with $N+\theta>2$, $l-\theta>-2$,  and $\Omega$ is a bounded or
unbounded domain. Through a suitable transformation of the form
$v(x)=|x|^\sigma u(x)$, equation (1) can be rewritten as a
 nonlinear Schr\"odinger equation with Hardy
potential
$$-\Delta u=|x|^\alpha |u|^{p-1}u+\frac{\ell}{|x|^2} u \;\; \mbox{in $\Omega \subset \R^N \;\; (N \geq 2)$}, \leqno{(2)}$$
where $p>1$, $\alpha \in (-\infty, \infty)$ and $\ell \in (-\infty,
(N-2)^2/4)$.

We show that under our chosen setting for the finite Morse index
theory of (1), the stability of a solution to (1) is unchanged under
various natural transformations. This enables us to reveal two
critical values of the exponent $p$ in (1) that divide the behavior
of finite Morse index solutions of (1), which in turn yields two
critical powers for (2)  through the transformation. The latter
appear difficult to obtain by working directly with (2).

\end{abstract}

\maketitle \baselineskip 18pt
\section{Introduction and main results}
\setcounter{equation}{0} In this paper, we study the properties of
finite Morse index solutions to the following weighted nonlinear
elliptic equation
\begin{equation}
\label{p*} -\mbox{div} (|x|^\theta \nabla v)=|x|^l |v|^{p-1}v \;\;\;
\mbox{in $\Omega \subset \R^N \; (N \geq 2)$},
\end{equation} where
$p>1$, $\theta, l\in\R^1$,  and $\Omega$ is a bounded or
unbounded domain. We are particularly interested in the cases that $\Omega$ is a punctured ball $B_R (0)
\backslash \{0\}$, an exterior domain $\R^N \backslash B_R$, or the
entire space $\R^N$. Here and throughout this paper, we use $B_r
(x)$ to denote the open ball in $\R^N$ centered at $x$ with radius
$r$. We also write $B_r=B_r (0)$.

An interesting classification of finite Morse index solutions to this equation in the case
$\Omega=\R^N$ (or $\R^N\setminus\{0\}$) and $\theta=l=0$ was given by Farina \cite{Fa} recently.
More recently such solutions in the case $\theta=0$ and $l>-2$ were considered in \cite{DDG, WY}.
Other recent related research on finite Morse index solutions can be found in
\cite{Da1, Da2, Da4, DY, DG, DF1, DF2}, where more references  are given.

Our interest in the  general case of \eqref{p*} was partly motivated
by research on the following
 nonlinear Schr\"odinger equation with Hardy
potential,
$$ -\Delta u=|x|^\alpha |u|^{p-1}u+\frac{\ell}{|x|^2} u \;\; \mbox{in $\Omega \subset \R^N$ ($N \geq
2$)},  \leqno(P)$$ where $p>1$, $\alpha \in (-\infty, \infty)$,
$\ell \in (-\infty, (N-2)^2/4)$. Equations of this type (with $N\geq
3$) arise in the study of nonlinear Schr\"odinger equations when the
field presents a (possible) singularity at the origin and have
attracted extensive studies in the past three decades; see, for
example, \cite{VV, DdPMW, FT, RW, Sm, SW, Te} and the references
therein.

 If we define
 \begin{equation}
\label{2'} v(x)=|x|^\sigma u(x), \;\;\;\;
\sigma=\frac{N-2}{2}-\sqrt{\Big(\frac{N-2}{2} \Big)^2-\ell},
\end{equation}
then (P) is reduced to \eqref{p*} with $\theta=-2 \sigma$ and
$l=\alpha-\sigma (p+1)$.
 Thus (P) can be reduced to \eqref{p*}, and vice versa.

 It should be noted that when $\theta\not=0$, the term
 $|x|^\theta$ in \eqref{p*} gives rise to a singularity (or
 degeneracy) at $x=0$ for the elliptic operator ${\rm div} (|x|^\theta \nabla v)$,
 and the notion of Morse index
 for solutions of
 \eqref{p*} need to be formulated with great care in order to make
 it consistent and useful. All the previous work on finite Morse index solutions
 that we are aware of dealt with elliptic equations with a uniformly elliptic operator.
 Therefore one might think that the form (P) is more
 natural to use than its equivalent equation in the form of
 \eqref{p*}. Our investigation here, however,
 suggests the opposite.

In this paper, we define the finite Morse index for \eqref{p*} in an
appropriate setting such that the stability of a solution to
\eqref{p*} is unchanged under several natural transformations. This
allows us to refine the calculations in \cite{DDG} to reveal two
critical values of $p$ for \eqref{p*} that divide the behavior of
finite Morse index solutions to \eqref{p*}, and through the
transformation \eqref{2'}, we obtain two critical powers for (P). As
will become clear below, the critical powers for (P) can be
expressed by relatively simple formulas in parameters appearing in
its equivalent form \eqref{p*}, but the formulas become very
complicated in terms of the parameters of (P) itself, which makes
them difficult to obtain by working on (P) directly. In a recent
work \cite{JL}, the methods of \cite{Fa} and \cite{DDG} are further
developed and applied to (P), but the calculations turn out to be
tedious in terms of the parameters appearing in (P), and the authors
have not found the optimal critical power for $p$ in the general
case.

To motivate some of the questions we investigate, and to get a taste
of how \eqref{p*} may be more natural to work with than (P), we
first recall two classical results of Bidaut-V\'eron and V\'eron
\cite{VV} (see Theorems 3.2, 3.3, 3.4 and Remark 3.2 in \cite{VV}).

\smallskip

\noindent
{\bf Theorem A.}
{\it
 Assume that $p \in (1, \infty) \backslash \{\frac{N+2+2
\alpha}{N-2} \}$,  $\ell<\frac{2+\alpha}{p-1}
\Big(N-2-\frac{(2+\alpha)}{p-1} \Big)$ and $u$ is a positive
solution of $(P)$ (with $N\geq 3$) in $B_R \backslash \{0\}$ such
that for some positive constant $C$,
$$|x|^{\frac{2+\alpha}{p-1}} u(x) \leq C \;\;\; \mbox{in $B_R \backslash \{0\}$}.$$
Then we have the following:

(i) either there exists $\eta>0$ such that
$$\lim_{x \to 0} u(x)|x|^{\frac{(N-2-\sqrt{(N-2)^2-4
\ell})}{2}}=\eta,$$

(ii) or there exists a positive solution $\omega$ of
$$\Delta_{S^{N-1}} \omega-\Big[\frac{2+\alpha}{p-1}
\Big(N-2-\frac{2+\alpha}{p-1}\Big)-\ell \Big]  \omega+\omega^p=0$$
on $S^{N-1}$ such that
$$\lim_{r \to 0} u(r\zeta)\,r^{\frac{2+\alpha}{p-1}} =\omega
(\zeta)$$ in the $C^k (S^{N-1})$-topology for any $k \in
\mathbb{N}$.
}

\smallskip

\noindent {\bf Theorem B.} {\it Assume that $p \in (1, \infty)
\backslash \{\frac{N+2+2 \alpha}{N-2} \}$,
$\ell<\frac{2+\alpha}{p-1} \Big(N-2-\frac{2+\alpha}{p-1} \Big)$ and
$u$ is a positive solution of $(P)$ (with $N\geq 3$) in $\R^N
\backslash B_R$ such that for some positive constant $C$,
$$|x|^{\frac{2+\alpha}{p-1}} u(x) \leq C \;\;\; \mbox{in $\R^N \backslash B_R$}.$$
Then we have the following:

(i) either there exists $\eta>0$ such that
$$\lim_{|x| \to \infty} u(x)|x|^{\frac{(N-2+\sqrt{(N-2)^2-4
\ell})}{2}}=\eta,$$

(ii) or there exists a positive solution $\omega$ of
$$\Delta_{S^{N-1}} \omega-\Big[\frac{2+\alpha}{p-1}
\Big(N-2-\frac{2+\alpha}{p-1}\Big)-\ell \Big]  \omega+\omega^p=0$$
on $S^{N-1}$ such that
$$\lim_{r \to \infty} u(r\zeta)\, r^{\frac{2+\alpha}{p-1}} =\omega
(\zeta)$$ in the $C^k (S^{N-1})$-topology for any $k \in
\mathbb{N}$.
}
\smallskip

If $1<p<\frac{N+2}{N-2}$, then the estimate
$|x|^{\frac{2+\alpha}{p-1}} u(x) \leq C$ in Theorems A and
B automatically holds for arbitrary $\alpha$ and $\ell$; see Theorem 6.3 in
\cite{VV} (for the special case $\ell=0$, this was first proved in \cite{GS}). The proof
for this fact is based on some useful integral estimates obtained from the
Bochner-Lichnerowicz-Weitzenb\"ock formula in $\R^N$.

For the case $\ell=0$, it was shown in
 \cite{DDG} that such estimate continues to hold for a larger range of $p$  provided that the solution has finite Morse index.
It would be interesting to see what happens for $\ell\not=0$. This
question will be answered as a consequence of some general results
in this paper for \eqref{p*}.

Let $v$ be a positive solution of \eqref{p*}.
If we define $r=|x|$, $\zeta=\frac{x}{|x|}$ and
$$
z(t, \zeta)=r^{\frac{2+l-\theta}{p-1}} v(r\zeta), \;\;\; t=\ln
r,$$ then $z(t, \zeta)$ satisfies the equation
$$
\begin{aligned}
z_{tt}+&\Big(N+\theta-2-\frac{2(2+l-\theta)}{p-1} \Big)
z_t+\Delta_{S^{N-1}} z \\
&-\frac{2+l-\theta}{p-1}
\Big[N+\theta-2-\frac{2+l-\theta}{p-1} \Big] z+|z|^{p-1} z=0.
\end{aligned}$$
One easily sees that the arguments in the proof of Theorems 3.2 and 3.3
in \cite{VV} still work for the above equation provided that
$$\frac{2+l-\theta}{p-1} \Big[N+\theta-2-\frac{2+l-\theta}{p-1}
\Big]>0,
$$
which is satisfied if
\begin{equation}
\label{N'-tau-p}
N+\theta>2,\; l-\theta>-2 \mbox{
and $p>\frac{N+l}{N+\theta-2}$.}
\end{equation}

Therefore, the proof of Theorems 3.2 and 3.3 in \cite{VV} yields the
following  result for \eqref{p*} (note that  $N=2$ is allowed here).

\begin{thm}
\label{C} Assume that \eqref{N'-tau-p} holds,  $p \neq
\frac{N+2+2l-\theta}{N+\theta-2}$,  and
$v$ is a positive solution of \eqref{p*} in $B_R \backslash \{0\}$ such
that for some positive constant $C$,
$$|x|^{\frac{2+l-\theta}{p-1}} v(x) \leq C \;\;\; \mbox{in $B_R \backslash \{0\}$}.$$
Then either $x=0$ is a removable singularity or it is a nonremovable
singularity and
$$r^{\frac{2+l-\theta}{p-1}} v(r\zeta) \to \varpi (\zeta) \;\; \mbox{as $r
\to 0$ unformly in $\zeta \in S^{N-1}$},$$ where
$\varpi$ is a positive solution of
\begin{equation}
\label{4} \Delta_{S^{N-1}} \varpi-\Big(\frac{2+l-\theta}{p-1} \Big)
\Big[N+\theta-2-\Big(\frac{2+l-\theta}{p-1} \Big) \Big]
\varpi+\varpi^p=0 \;\;\; \mbox{on $S^{N-1}$}.
\end{equation}
\end{thm}

\begin{thm}
\label{D} Assume that  \eqref{N'-tau-p} holds,  $p \neq
\frac{N+2+2l-\theta}{N+\theta-2}$, and
$v$ is a positive solution of \eqref{p*} in $\R^N \backslash B_R$ such
that for some positive constant $C$,
$$|x|^{\frac{2+l-\theta}{p-1}} v(x) \leq C \;\;\; \mbox{in $\R^N \backslash B_R$}.$$
Then either
$$|x|^{N-2+\theta} v(x) \to \gamma \;\;\; \mbox{as $|x| \to \infty$ for some
$\gamma>0$}$$ or
$$r^{\frac{2+l-\theta}{p-1}} v(r\zeta) \to \varpi (\zeta) \;\; \mbox{as $r
\to \infty$ unformly in $\zeta \in S^{N-1}$},$$ where
$\varpi (\omega)$ is a positive solution of \eqref{4}.
\end{thm}

\begin{rem}
\label{P-p*}
It is easily checked that the condition in Theorems A and B on $\ell$, namely
\begin{equation}
\label{alpha-ell}
\ell<\frac{2+\alpha}{p-1}\left(N-2-\frac{2+\alpha}{p-1}\right),
\end{equation}
is equivalent to   \eqref{N'-tau-p} with $\theta=-2\sigma$ and $l=\alpha-\sigma(p+1)$.
\end{rem}

\bigskip

We now introduce the setting in which the finite Morse index theory
for \eqref{p*} will be developed. This is a crucial first step for
the analysis of this paper. As mentioned before, we need to choose
the setting with great care in order to make the notion of finite
Morse index useful. In particular, we want the stability of a
solution to \eqref{p*} unchanged under various natural
transformations, including \eqref{2'}, the Kelvin transformation
\eqref{6} and the transformation \eqref{dual} given below.

 For $\theta\in\R^1$, we denote by $H^{1,\theta} (\Omega)$ the
space of functions $\varphi $ such that
\[
|x|^{\frac{\theta}{2}} \varphi\in L^2(\Omega), \;
|x|^{\frac{\theta}{2}} |\nabla\varphi|\in L^2(\Omega),
\]
with
  norm
\[
  \|\varphi \|=\left(\int_\Omega |x|^{\theta} (\varphi^2+ |\nabla \varphi|^2)dx\right)^{1/2}.
  \]
   $H^{1,\theta}_{loc}(\Omega)$ is defined in the obvious way, and
we use $H_c^{1,\theta}(\Omega)$ to denote the subspace of functions
in $H^{1,\theta}(\Omega)$ which have compact supports in $\Omega$.
If $0\not\in\Omega$, clearly
$H^{1,\theta}_{loc}(\Omega)=H^1_{loc}(\Omega)$ and
$H^{1,\theta}_c(\Omega)=H^1_c(\Omega)$. If further
$0\not\in\overline\Omega$ and $\Omega$ is bounded, then
$H^{1,\theta}(\Omega)=H^1(\Omega)$.

\begin{rem}
\label{space} \begin{itemize} \item[(i)] If $N\geq 3$ and
$\Omega\subset \R^N$ is bounded, then $u\in H^{1,\theta}(\Omega)$ if
and only if $|x|^{\frac{\theta}{2}}u\in H^1(\Omega)$. This is a
direct consequence of the identity
\[
\nabla(|x|^{\frac{\theta}{2}}u)=|x|^{\frac{\theta}{2}}\nabla
u+\frac{\theta}{2}\frac{ x}{|x|^2} \,|x|^{\frac{\theta}{2}}u
\]
and the fact that $|x|^{-1}\in L^2(\Omega)$ for $N\geq 3$ when
$\Omega$ is bounded.

\item[(ii)] If $N=2$ and $0\in\Omega$, then $|x|^{-1}\not\in L^2(\Omega)$,
and the above statement is not true.
\item[(iii)]
If $N\geq 2$, $u\in L_{loc}^\infty(\Omega)$ and $N+\theta>2$, then
$u\in H_{loc}^{1,\theta}(\Omega)$ if and only if
$|x|^{\frac{\theta}{2}}u\in H_{loc}^{1}(\Omega)$. To see this, it
suffices to check that $|x|^{\frac{\theta}{2}-1}u\in
L_{loc}^2(\Omega)$ under the given conditions. Indeed, from
$N+\theta>2$ one obtains $|x|^{\frac{\theta}{2}-1}\in
L_{loc}^2(\Omega)$, which implies $|x|^{\frac{\theta}{2}-1}u\in
L_{loc}^2(\Omega)$ since by assumption $u\in
L_{loc}^\infty(\Omega)$.
\end{itemize}
\end{rem}

 We say that $v$ is a solution of \eqref{p*} if $v \in
H^{1, \theta}_{loc} (\Omega)\cap L^\infty_{loc}(\Omega)$ and
\begin{equation}
\label{5} \int_\Omega \Big(|x|^\theta \nabla v \cdot \nabla
\phi-|x|^l |v|^{p-1}v \phi\Big)=0 \;\;\;\; \forall \phi \in
H^{1,\theta}_c(\Omega)\cap L_{loc}^\infty(\Omega).
\end{equation}
Let us observe that if $v$ is a solution of \eqref{p*}, then by
standard elliptic regularity $v \in C^2(\Omega \backslash \{0\})$
and hence is a classical solution of \eqref{p*} in $\Omega
\backslash \{0\}$. In particular, $v \in C^2 (\Omega)$ whenever $0
\not \in \Omega$.  If $0\in\Omega$, then \eqref{5} has a hidden
restriction on $v$ at $x=0$ since $\int_\Omega |x|^l|v|^{p-1}v
\phi\, dx$ need not be defined for arbitrary $v\in
H_{loc}^{1,\theta}(\Omega)\cap L_{loc}^\infty(\Omega)$ and $\phi\in
H_c^{1,\theta}(\Omega)\cap L_{loc}^\infty(\Omega)$. However, this
hidden restriction disappears when $N+\theta>2$ and $l-\theta>-2$,
since in such a case, $N+l>0$ and $|x|^l\in L_{loc}^1(\Omega)$.

A solution $v$ of \eqref{p*} is said to be {\it stable} if
$$Q_v (\psi):=\int_\Omega \Big(|x|^\theta |\nabla \psi|^2-p |x|^l
|v|^{p-1} \psi^2\Big) \geq 0 \;\;\;\; \forall \psi \in
H^{1,\theta}_c (\Omega)\cap L_{loc}^\infty(\Omega).$$

Similar to \cite{Da1}, we say a solution of \eqref{p*} has {\it
Morse index} $k \geq 0$ if $k$ is the maximal dimension of all
subspaces $X_k$ of $H_c^{1,\theta}(\Omega)\cap
L_{loc}^\infty(\Omega)$ such that $Q_v(\psi)<0$ for any $\psi\in
X_k\setminus\{0\}$. Thus $v$ is stable if and only if it has Morse
index 0. Moreover, if  $v$ has finite Morse index over the domain
$\Omega$, then there exists a compact subset $\mathcal{K}$ of
$\Omega$ such that $v$ is stable over any domain
$\Omega'\subset\Omega\setminus \mathcal{K}$.

The above setting allows us to establish the following integral
estimate for stable solutions of \eqref{p*}, which is a key step for
the success of this approach. This estimate is an extension of  Proposition 4
in \cite{Fa} (for $\theta=l=0$) and  Proposition 1.7 in \cite{DDG}
(for $\theta=0$ and $l>-2$), albeit that we  have added an extra term $|\psi|\frac{|\nabla\psi|}{|x|}$ in the right hand side.
However, this extra term does not affect the key estimates in its applications, even for the special cases considered in \cite{DDG} and \cite{Fa}.

\begin{prop}
\label{K} Let $\Omega$ be a domain (bounded or not) of $\R^N \; (N
\geq 2)$. Let $v \in H_{loc}^{1,\theta} (\Omega)\cap L^\infty_{loc}(\Omega)$ be a stable solution
of \eqref{p*} with $p>1$. Then for any $\gamma \in [1, 2p+2
\sqrt{p(p-1)}-1)$ and any integer $m \geq \max
\{\frac{p+\gamma}{p-1},2\}$ there exists a constant $C>0$ depending
only on $p,m, \gamma,l$ and $\theta$ such that
\begin{equation}
\begin{aligned}
\label{10} &\int_\Omega \Big(|x|^\theta |\nabla
(|v|^\frac{\gamma-1}{2} v)|^2+|x|^l |v|^{\gamma+p} \Big) |\psi|^{2m}\\
&\leq C \int_\Omega |x|^{\frac{\theta (\gamma+p)-l (\gamma+1)}{p-1}}
\Big(|\nabla \psi|^2+|\psi||\Delta \psi|+|\psi| \frac{|\nabla
\psi|}{|x|} \Big)^{\frac{p+\gamma}{p-1}}
\end{aligned}
\end{equation}
for all test functions $\psi \in C_0^2 (\Omega)$ satisfying $|\psi|
\leq 1$ in $\Omega$.
\end{prop}

As in \cite{DDG}, the Kelvin transformation will be a useful tool in this paper.
If $v$ is a solution of \eqref{p*} over $B_R \backslash \{0\} \; (N
\geq 2)$, then the function $w$ defined by the Kelvin
transformation
\begin{equation}
\label{6} w(y)=|x|^{N-2+\theta} v (x), \;\;\;\;\;\;
y=\frac{x}{|x|^2}
\end{equation}
satisfies the equation
\begin{equation}
\label{7} -\mbox{div} (|y|^\theta \nabla w)=|y|^\beta |w|^{p-1} w
\;\;\;\; \mbox{for $y \in \R^N \backslash {\overline {B_{1/R}}}$},
\end{equation}
with $\beta=(N-2+\theta)(p-1)-(4+l-2 \theta)$. We notice that
$\beta-\theta>-2$ when $p>\frac{N+l}{N-2+\theta}$.

We have the following proposition which shows that the Kelvin
transformation in \eqref{6} keeps the stability of $v$.

\begin{prop}
\label{F} A solution $v$ of \eqref{p*} is stable in $B_R \backslash
\{0\}$ if and only if the function $w$ obtained by the Kelvin
transformation in \eqref{6} is a stable solution of \eqref{7} in
$\R^N \backslash {\overline {B_{1/R}}}$.
\end{prop}

The next proposition discusses the stability property between solutions of (P)
and \eqref{p*}. Recall that these two problems are related through
$v(x)=|x|^\sigma u(x)$ with $\theta=-2 \sigma$,
$l=\alpha-\sigma (p+1)$ and $\sigma=\frac{1}{2} [N-2-\sqrt{(N-2)^2-4
\ell}]$.

We say $u$ is a solution of (P) if $u\in H^1_{loc}(\Omega)$,
$|x|^\sigma u\in L^\infty_{loc}(\Omega)$, and
\[
\int_\Omega \Big(\nabla u\cdot\nabla \phi-\ell |x|^{-2}u
\phi-|x|^\alpha |u|^{p-1}u \phi\Big)=0 \;\;\;\forall \phi\in
H^1_c(\Omega) \mbox{ with } |x|^\sigma\phi\in
L_{loc}^\infty(\Omega).
\]
It is said to be {\it stable} if
$${\mathcal Q}_u (\phi):=\int_\Omega \Big(|\nabla \phi|^2-\ell
|x|^{-2} \phi^2-p |x|^\alpha |u|^{p-1} \phi^2\Big) \geq 0 \ $$ for
all $\phi \in H_c^1 (\Omega) \mbox{ with } |x|^\sigma\phi\in
L_{loc}^\infty(\Omega).$

We say that $u$ has Morse index $k\geq 0$ if $k$ is the maximal
dimension of all subspaces $Y_k$ of $Y:=\{\phi\in H_c^1(\Omega):
|x|^\sigma \phi\in L_{loc}^\infty(\Omega)\}$ such that ${\mathcal
Q}_u (\phi)<0$ for any $\phi\in Y_k\setminus \{0\}$.

\begin{prop}
\label{p5.2} Let $u$ be a solution of
$(P)$. Then $v(x):=|x|^\sigma u(x) $
is a stable solution of \eqref{p*} if and only if $u$ is a stable
solution of $(P)$.
\end{prop}

To introduce the other results of this paper, we need to define two
critical powers for \eqref{p*}. In order to use calculations  in
\cite{DDG} by similarity, in the following, we denote
$$N'=N+\theta \;\;\; \mbox{and} \;\;\;\; \tau=l-\theta$$
for fixed $l$ and $\theta$ in $\R^1$. We assume from now on that
\begin{equation}
\label{7-1} N'>2 \;\;\; \mbox{and} \;\;\; \tau>-2,
\end{equation}
 unless otherwise specified.

 To better understand the above restriction on $N'$ and $\tau$, we make use of another transformation
 \begin{equation}
 \label{dual}
 z(y)=v(x),\; y=\frac{x}{|x|^2}.
 \end{equation}
 A simple calculation shows that under this transformation $v$ is a solution to \eqref{p*} if and only if $z$
 is a solution to
 \begin{equation}
 \label{dual-p*}
 -\mbox{div}(|y|^{\tilde \theta}\nabla z)=|y|^{\tilde l} |z|^{p-1}z,\;\; |y|^{-2}y\in\Omega,
 \end{equation}
 with
 \[
 \tilde\theta=4-2N-\theta,\; \tilde l=-2N-l.
 \]
If we define
\[
\tilde N':=N+\tilde\theta,\;\; \tilde \tau:=\tilde l-\tilde\theta,
\]
then $\tilde N'+N'=4$ and $\tilde\tau+\tau=-4$. Thus
\[
 \mbox{$N'<2$ \ \  if and only if\ \  $\tilde N'>2$},\]  and
\[ \mbox{
 $\tau<-2$\ \  if and only if\ \  $\tilde \tau>-2$}.
\]
Moreover, the stability of the solution of \eqref{p*} is unchanged
under the transformation \eqref{dual}:

\begin{prop} \label{L} A solution $v$ of \eqref{p*} is stable in $B_R \backslash
\{0\}$ if and only if the function $z$ obtained by the
transformation in \eqref{dual} is a stable solution of
\eqref{dual-p*} in $\R^N \backslash {\overline {B_{1/R}}}$.
\end{prop}

Thus for every  result we obtain in the case of \eqref{7-1} there is
a parallel result  in the case of $N'<2$ and $\tau<-2$ through the
transformation \eqref{dual}.

We will show that if $N'\geq 2$ and $\tau\leq -2$, then  problem
\eqref{p*} has no positive solution over any punctured ball $B_R
\backslash \{0\}$ (see Theorem 4.2 below). This implies, by the
Kelvin transformation, problem \eqref{p*} has no positive solution
over any exterior domain $\R^N \backslash B_R$ if  $p \leq
\frac{N'+\tau}{N'-2}$. This also implies, by the transformation
\eqref{dual}, that problem \eqref{p*} has no positive solution over
any exterior domain $\R^N \backslash B_R$ if $N'\leq 2$ and $\tau\geq
-2$.

For these reasons, the case $N'\geq 2$ and $\tau\leq -2$, and the case
$N'\leq 2$ and $\tau\geq -2$, are not considered further.\footnote{Note,
however, for $N'\geq 2$ and $\tau\leq -2$, one may still consider
\eqref{p*} over an exterior domain, and for $N'\leq 2$ and $\tau\geq
-2$, one may consider \eqref{p*} over a punctured ball. But we will
not pursue these cases here.} Our focus will be mainly on the case
\eqref{7-1}.

Suppose \eqref{7-1} holds and let
$$f(p)=p \frac{2+\tau}{p-1}
\Big(N'-2-\frac{2+\tau}{p-1} \Big).$$ Evidently,
$$f \Big(\frac{N'+\tau}{N'-2} \Big)=0, \;\;\;\;
f(\infty)=(2+\tau)(N'-2).$$

Replacing $(N,\alpha)$ in the calculations on page 3285 of
\cite{DDG} by $(N',\tau)$, we find that the equation
$f(p)=\frac{(N'-2)^2}{4}$ always has a solution in the interval
$(\frac{N'+\tau}{N'-2}, \frac{N'+2+2\tau}{N'-2})$. We denote this
solution by $P_-(N',\tau)$. A simple calculation shows that
$f(p)=\frac{(N'-2)^2}{4}$ is equivalent to
\[ ap^2-bp+cp=0\]
with
\begin{equation}
\label{abc}
\left\{\begin{array}{l}
a=(N'-2)(N'-4\tau-10),
\vspace{0.2cm}\\
 b=2(N'-2)^2-4(\tau+2)(\tau+N'),
 \vspace{0.2cm}\\
c=(N'-2)^2.
\end{array}
\right.
\end{equation}
From this, we obtain
$$P_- (N', \tau):=\frac{(N'-2)^2-2(2+\tau)(N'+\tau)-2(2+\tau)
\sqrt{(2+\tau)(2N'+\tau-2)}}{(N'-2)(N'-4 \tau-10)}$$
if $N'\not=4\tau+10$, and $P_-(N',\tau)=\frac{4}{3}$ if $N'=4\tau+10$.
Moreover, when
$2<N' \leq 10+4 \tau$, $f$ has the property
\begin{equation}
\label{8} \left \{ \begin{array}{ll} f(p)<\frac{(N'-2)^2}{4}
\;\;\;\; &\mbox{for $1<p<P_-(N', \tau)$},\\
f(p)>\frac{(N'-2)^2}{4} \;\;\;\;&\mbox{for $p>P_-(N', \tau)$}.
\end{array} \right.
\end{equation}
When $N'>10+4\tau$, there exists a second root of $f(p)=\frac{(N'-2)^2}{4}$ in $(1,\infty)$, given by
$$P_+ (N', \tau):=\frac{(N'-2)^2-2
(2+\tau)(N'+\tau)+2(2+\tau) \sqrt{(2+\tau)(2N'+\tau-2)}}{(N'-2)(N'-4
\tau-10)},$$ and it has the properties
\[
\frac{N'+2+2\tau}{N'-2}<P_+(N',\tau)<\infty,
\]
and
\begin{equation}
\label{9} \left \{ \begin{array}{ll} f(p)<\frac{(N'-2)^2}{4}
\;\;\;\; &\mbox{for $p \in (1, P_- (N', \tau)) \cup
(P_+ (N', \tau), \infty)$},\\
f(p)>\frac{(N'-2)^2}{4} \;\;\;\;& \mbox{for $p \in (P_- (N', \tau),
P_+ (N', \tau))$}.
\end{array} \right.
\end{equation}

We will show that the number
$$p_c (N', \tau)=\left \{ \begin{array}{ll} \infty \;\;\;\;
&\mbox{if $2<N' \leq 10+4 \tau$},\\
P_+ (N', \tau) \;\;\;\; &\mbox{if $N'>10+4 \tau$},
\end{array} \right.$$
serves as a critical power for \eqref{p*}. The number
\[
\tilde p_c(N',\tau):=P_-(N',\tau)<\frac{N'+2+2\tau}{N'-2}<p_c(N',\tau)
\]
is also a critical value for \eqref{p*}.

The first important role played by $p_c(N', \tau)$ can be seen from the
following theorem.

\begin{thm}
\label{F'} If $2<N' \leq 10+4 \tau$ and $p>1$, or $N'>10+4 \tau$ and
$1<p<p_c (N',\tau)$, and if $v\in H^{1,\theta}_{loc}(\R^N)\cap L^\infty_{loc}(\R^N)$ is a
stable solution of \eqref{p*} (nonnegative or not) in $\R^N$ $(N\geq 2)$, then $v
\equiv 0$; on the other hand, if $p \geq p_c (N',\tau)$, \eqref{p*}
admits a family of stable positive radial solutions in $\R^N$.
\end{thm}

For the special case $\theta=l=0$, the above result was first
obtained in \cite{Fa}. When $\theta=0$ and $l>-2$, it was proved in
\cite{DDG}. See \cite{DLF} for the important role played by $p_c(N',\tau)$
on the behavior of radially symmetric solutions.

 All the other results in this paper treat
equations over a punctured domain or an exterior domain. We say that
a positive solution $v$ of \eqref{p*} has an {\it isolated
singularity} at 0 if $\Omega$ contains a punctured ball $B_r
\backslash \{0\}$, $0 \not \in \Omega$ and $v$ tends to $\infty$
along some sequence $x_n \to 0$. If on the other hand $\lim_{|x| \to
0} v(x)=\gamma$ is a finite number and $v$ becomes a positive
solution of \eqref{p*} over $B_r$ upon defining $v(0)=\gamma$, we
say that $x=0$ is a removable singularity of $v$.

Let $\Omega^* \subset \R^N \; (N \geq 3)$ be a bounded domain such
that $0 \in \Omega^*$. A positive solution $v$ of \eqref{p*} in
$\R^N \backslash \Omega^*$ is called a {\it fast decay} solution if
$\lim_{|x| \to \infty} |x|^{N+\theta-2} v(x)=\gamma$ for some
$\gamma>0$.

The following two results give sufficient conditions to meet the
requirements in Theorems \ref{C} and \ref{D}. We note that in these
two theorems, we have no restriction on $\theta$ and $l$.

\begin{thm}
\label{G} Let $\Omega_0 \subset \R^N \; (N \geq 2)$ be a bounded
domain containing 0, and let $v$ be a positive solution of \eqref{p*}
in $\Omega_0 \backslash \{0\}$ with arbitrary $\theta,\;l\in \R^1$. If $v$ has finite Morse index and
$1<p<p_c (N,0)$, then there exist $C>0$ and $\epsilon>0$ such that
$$|x|^{\frac{2+\tau}{p-1}} v(x) \leq C \;\;\; \mbox{for
$0<|x|<\epsilon$}.$$ Hence by Theorem \ref{C}, when $p \in
(\frac{N'+\tau}{N'-2},  p_c (N,0)) \backslash \{\frac{N'+2+2
\tau}{N'-2}\}$,
$$ \left \{ \begin{array}{l} \mbox{either $v$ has a removable
singularity at $x=0$, or}\\
\mbox{$r^{\frac{2+\tau}{p-1}} v(r\zeta) \to \varpi (\zeta)$ as $r
\to 0$ uniformly in $\zeta\in S^{N-1}$},
\end{array} \right. \leqno{(A_0)}$$
where $\varpi$ is a positive solution of \eqref{4}.
\end{thm}

\begin{thm}
\label{H} Let $\Omega_0 \subset \R^N\; (N \geq 2)$ be a bounded
domain containing 0, and let $v$ be a positive solution of \eqref{p*}
in $\R^N \backslash \Omega_0$ with arbitrary $\theta,\;l\in \R^1$. If $v$ has finite Morse index and
$1<p<p_c (N,0)$, then there exist $C>0$ and
$\epsilon>0$ such that
$$|x|^{\frac{2+\tau}{p-1}} v(x) \leq C \;\;\;\; \mbox{for
$|x|>\epsilon^{-1}$}.$$ Hence by Theorem \ref{D}, when $p \in
(\frac{N'+\tau}{N'-2}, p_c (N,0)) \backslash \{\frac{N'+2+2
\tau}{N'-2}\}$,
$$ \left \{ \begin{array}{l} \mbox{either $v$ is a fast decay
solution, i.e., $\lim_{|x|\to\infty} |x|^{N'-2}v(x)=\gamma>0$, or}\\
\mbox{$r^{\frac{2+\tau}{p-1}} v(r\zeta) \to \varpi (\zeta)$ as $r
\to \infty$ uniformly in $\zeta \in S^{N-1}$},
\end{array} \right. \leqno{(A_\infty)}$$
where $\varpi$ is a positive solution of \eqref{4}.
\end{thm}

For the special case $\theta=0$ and $l>-2$, the above two theorems
were first proved in \cite{DDG}.

With more restrictions on $p$, we can determine the alternatives in
$(A_0)$ and $(A_\infty)$.

\begin{thm}
\label{I} Let $\Omega_0 \subset \R^N \; (N \geq 2)$ be a bounded
domain containing 0, and let $v$ be a positive solution of \eqref{p*}
in $\Omega_0 \backslash \{0\}$. If $v$ has finite Morse index and if
\begin{equation}
\label{cond-p}
\tilde p_c(N',\tau)<p<\min\{p_c (N',\tau), p_c(N,0)\},\; p\not=\frac{N'+2+2\tau}{N'-2},
\end{equation}
then $x=0$ must be a removable singularity of $v$.

On the other hand, for $p \geq p_c (N',\tau)$ or $p\in
(\frac{N'+\tau}{N'-2},\tilde p_c(N',\tau))$, \eqref{p*} has a
positive stable solution on $\R^N \backslash \{0\}$ with an isolated
singularity at 0 (which is $V_\infty$ given below).
\end{thm}

\begin{rem} We will show in Remark \ref{r3.1} that
the function $p_c (\cdot, \tau)$ is a decreasing function for fixed
$\tau$ and $p_c (N', \cdot)$ is an increasing function for fixed
$N'$, as long as the value of the functions is finite (i.e.,
$N'>10+4\tau$). Moreover, when $\tau = \frac{p-1}{2p+2
\sqrt{p(p-1)}}\,\theta$, we have $p_c(N',\tau)=p_c(N,0)$. Therefore
\[
\min\left\{p_c (N',\tau), p_c(N,0)\right\}=\left\{
\begin{array}{ll}
p_c(N',\tau) & \mbox{ if } \tau \leq \frac{p-1}{2p+2 \sqrt{p(p-1)}}\,
\theta,\vspace{0.3cm}\\
p_c(N,0) & \mbox{ if } \tau > \frac{p-1}{2p+2 \sqrt{p(p-1)}}\,
\theta.
\end{array}\right.
\]
\end{rem}
Note that the inequality $\tau \leq \frac{p-1}{2p+2 \sqrt{p(p-1)}}\,
\theta$ is equivalent to
\begin{equation}
\label{9-1} l \leq \left(1+\frac{p-1}{2p+2 \sqrt{p(p-1)}}\right)
\theta.
\end{equation}
\begin{thm}
\label{J} Let $\Omega_0 \subset \R^N \; (N \geq 2)$ be a bounded
domain containing 0, and let $v$ be a positive solution of
\eqref{p*} in $\R^N \backslash \Omega_0$. If $v$ has finite Morse
index and if \eqref{cond-p} holds, then $v$ must be a fast decay
solution $(i.e., \lim_{|x|\to\infty} |x|^{N'-2}v(x)=\gamma>0)$.

On the other hand, for $p \geq p_c (N',\tau)$ or $p\in
(\frac{N'+\tau}{N'-2},\tilde p_c(N',\tau))$, \eqref{p*} has a stable
positive radial solution on $\R^N \backslash \{0\}$ which decays at
the slower rate $|x|^{-\frac{2+\tau}{p-1}}$ at $\infty$ (which is
$V_\infty$ given below).
\end{thm}

\begin{rem}
\label{correction}
Theorems \ref{I} and \ref{J} indicate that the conclusions in Theorems 1.5 and 1.6 of \cite{DDG}
hold only for $p$ in the range
\[
\underline{p}(\alpha)<p<\overline p(\alpha^-),\; p\not=\frac{N+2+2\alpha}{N-2}
\]
instead of
\[
\frac{N+\alpha}{N-2}<p<\overline p(\alpha^-),\; p\not=\frac{N+2+2\alpha}{N-2}
\]
as claimed there. The mistake in \cite{DDG} is caused by the statement that
\[
f\left(\frac{N+2+2\beta}{N-2}\right)>\frac{(N-2)^2}{4} \mbox{ implies } p>\underline p(\beta).
\]
The above statement is true if $\beta$ is independent of $p$, but
$\beta=(N-2)(p-1)-(4+\alpha)$ in \cite{DDG}. We also note that
$\underline p(\alpha)$ is increasing in $\alpha\in (-2,\infty)$
instead of decreasing (as stated in \cite{DDG}).
\end{rem}

It is easily checked  that
$$V_\infty (x)=C_0 |x|^{-\frac{2+\tau}{p-1}}, \;\;\;\;\mbox{ with } C_0=\Big
\{\frac{2+\tau}{p-1} \Big(N'-2-\frac{2+\tau}{p-1} \Big) \Big
\}^{1/(p-1)},$$ is a positive radial solution of \eqref{p*} over
$\R^N\setminus\{0\}$ provided that $\tau>-2$ and
$p>\frac{N'+\tau}{N'-2}$. For calculation convenience, we note that
if $v=v(r)$ is a radial solution of \eqref{p*}, then $v(r)$
satisfies
$$v_{rr}+\frac{N'-1}{r} v_r+r^\tau |v|^{p-1} v=0.$$
Moreover, we will show that $V_\infty$ is the only positive radial solution of \eqref{p*} over a
punctured ball $B_R \backslash \{0\}$ that has a singularity at 0 if
$p>\frac{N'+2+2 \tau}{N'-2}$ (see Theorem 4.3 below).

\begin{rem} Theorems \ref{I} and \ref{J} imply that, if $\tilde
p_c(N',\tau)<p<p_c (N', \tau)$, $p\not=\frac{N'+2+2\tau}{N'-2}$ and
$\tau\leq \frac{p-1}{2p+\sqrt{p(p-1)}}\,\theta$, then the Morse
index of $V_\infty$ is $\infty$ as a positive solution of \eqref{p*}
over any punctured ball $B_r \backslash \{0\}$, or over any $\R^N
\backslash B_R$, but when $p \geq p_c (N', \tau)$ or $p\in
(\frac{N'+\tau}{N'-2},\tilde p_c(N',\tau)]$, the Morse index of
$V_\infty$ is reduced to 0. We do not know whether  Theorems \ref{I}
and \ref{J} still hold if $\min\{p_c(N',\tau),p_c(N,0)\}$ in
\eqref{cond-p} is replaced by $p_c(N',\tau)$ when
$\tau>\frac{p-1}{2p+\sqrt{p(p-1)}}\,\theta$.
\end{rem}

The rest of the paper is organized in the following way. In Section
2, we give the proofs of Propositions \ref{K}, \ref{F}, \ref{p5.2}
and \ref{L}. In Section 3, we prove Theorem \ref{F'}. Section 4 is
devoted to the proof of Theorems \ref{G} and \ref{H}, while Theorems
\ref{I} and \ref{J} are proved in Section 5, the last section of the
paper.

\section{Proofs of the basic results}
\setcounter{equation}{0}

In this section, we collect the proofs of all the basic results which will serve as
tools  in the proofs of our other results.

\subsection{Proof of Proposition \ref{K}.} We follow the lines  of
the proof of Proposition 4 in \cite{Fa} and Proposition 1.7 of
\cite{DDG}, but with considerable modifications. We divide the proof into three steps.

{\it Step 1.} For any $\varphi \in C_0^2 (\Omega)$,
\begin{equation}
\label{5.1} \int_\Omega |x|^\theta |\nabla (|v|^{\frac{\gamma-1}{2}}
v)|^2 \varphi^2=\frac{(\gamma+1)^2}{4 \gamma} \int_\Omega |x|^l
|v|^{p+\gamma} \varphi^2+\frac{\gamma+1}{4 \gamma} \int_\Omega
|v|^{\gamma+1} \mbox{div} (|x|^\theta \nabla (\varphi^2)).
\end{equation}
This is obtained by taking $\phi=|v|^{\gamma-1}v \varphi^2$ in
\eqref{5}.

{\it Step 2.} For any $\varphi \in C_0^2 (\Omega)$, we have
\begin{equation}
\label{5.2}
\begin{aligned} & \Big(p-\frac{(\gamma+1)^2}{4 \gamma} \Big) \int_\Omega
|x|^l |v|^{\gamma+p} \varphi^2\\
& \leq \int_\Omega |x|^\theta
|v|^{\gamma+1} |\nabla \varphi|^2+\frac{\gamma-1}{4 \gamma}
\int_\Omega |v|^{\gamma+1} \mbox{div} (|x|^\theta \nabla
(\varphi^2)).
\end{aligned}
\end{equation}

The function $\psi=|v|^{\frac{\gamma-1}{2}} v \varphi$ belongs to
$H_{c}^{1,\theta}(\Omega)\cap L_{loc}^\infty(\Omega)$, thus it can
be used as a test function in the quadratic inequality $Q_v
(\psi)\geq 0$. Taking this test function and using \eqref{5.1}, we
can easily obtain \eqref{5.2}.

{\it Step 3.} For any $\gamma \in [1, 2p+2 \sqrt{p(p-1)}-1)$ and any
$m \geq \max \{\frac{p+\gamma}{p-1},2\}$, there exists a constant
$C>0$ depending only on $p,m,\gamma,l, \theta$ such that
\begin{equation}
\label{5.3} \int_\Omega |x|^l |v|^{p+\gamma} |\psi|^{2m} \leq C
\int_\Omega |x|^{\frac{\theta (\gamma+p)-l(\gamma+1)}{p-1}}
\Big(|\nabla \psi|^2+|\psi||\Delta \psi|+|\psi|\frac{|\nabla
\psi|}{|x|} \Big)^{\frac{p+\gamma}{p-1}}
\end{equation}
for all test function $\psi \in C_0^2 (\Omega)$ satisfying $|\psi|
\leq 1$ in $\Omega$.

From \eqref{5.2} we see that for any $\varphi \in C_0^2 (\Omega)$,
\begin{equation}
\label{5.4} \eta \int_\Omega |x|^l |v|^{p+\gamma} \varphi^2 \leq
\kappa \int_\Omega |v|^{\gamma+1} \mbox{div} (|x|^\theta \nabla
(\varphi^2))+\int_\Omega |x|^\theta |v|^{\gamma+1} |\nabla
\varphi|^2
\end{equation}
with
$$\eta=p-\frac{(\gamma+1)^2}{4 \gamma}, \;\;\;\;
\kappa=\frac{\gamma-1}{4 \gamma}.$$

For any $\gamma \in [1, 2p+2 \sqrt{p(p-1)}-1)$, an elementary
analysis shows that $\eta>0$.

For any $\psi \in C_0^2 (\Omega)$ with $|\psi| \leq 1$ in $\Omega$,
we set $\varphi=\psi^m$. Since $m \geq 2$, the function $\varphi$
belongs to $C_0^2 (\Omega)$ and it follows from \eqref{5.2} that
$$ I:=\int_\Omega |x|^l |v|^{\gamma+p} |\psi|^{2m} \leq C_{m,p,\gamma,l,
\theta}^1 \int_\Omega |x|^\theta |v|^{\gamma+1} |\psi|^{2m-2}
\Big(|\nabla \psi|^2+|\psi||\Delta \psi|+|\psi|\frac{|\nabla
\psi|}{|x|} \Big)$$
 where $C_{m,p,\gamma,l,\theta}^1>0$ depends on
$m,p,\gamma,l,\theta$. An application of H\"older's inequality
yields
\[
\begin{aligned}
&\int_\Omega |x|^\theta |v|^{\gamma+1} |\psi|^{2m-2} \Big(|\nabla
\psi|^2+|\psi||\Delta \psi|+|\psi|\frac{|\nabla \psi|}{|x|} \Big)\\
&\;\leq I^{\frac{\gamma+1}{\gamma+p}} \Big[\int_\Omega |x|^{\frac{\theta
(\gamma+p)-l (\gamma+1)}{p-1}} |\psi|^{2(m-\frac{\gamma+p}{p-1})}
\Big (|\nabla \psi|^2+|\psi||\Delta \psi|+|\psi| \frac{|\nabla
\psi|}{|x|} \Big)^{\frac{p+\gamma}{p-1}} \Big
]^{\frac{p-1}{\gamma+p}}
\end{aligned}
\]
and hence
\begin{equation}
\label{5.5} \int_\Omega |x|^\theta |v|^{\gamma+p} |\psi|^{2m} \leq C
\int_\Omega |x|^{\frac{\theta (\gamma+p)-l (\gamma+1)}{p-1}} \Big
(|\nabla \psi|^2+|\psi||\Delta \psi|+|\psi| \frac{|\nabla
\psi|}{|x|} \Big)^{\frac{p+\gamma}{p-1}},
\end{equation}
which proves \eqref{10} and Proposition \ref{K}. \qed

\subsection{Proof of Proposition \ref{F}} For any given $\psi \in H^1_c (B_R \backslash \{0\})\cap L_{loc}^\infty(B_R \backslash \{0\})$,
we define
$${\tilde \psi} (y)=|x|^{N'-2} \psi (x), \;\;\;\;y=\frac{x}{|x|^2}.$$
Clearly ${\tilde \psi} \in H^1_c (\R^N \backslash {\overline
{B_{1/R}}})\cap L_{loc}^\infty(\R^N \backslash {\overline
{B_{1/R}}})$. (Recall that for this kind of domains
$H^{1,\theta}_c(\Omega)=H^1_c(\Omega)$.) Moreover,
\begin{eqnarray*}
& &\int_{\R^N \backslash {\overline {B_{1/R}}}} \Big[|y|^\theta
|\nabla_y {\tilde \psi}|^2-p |y|^\beta |w|^{p-1} {\tilde \psi}^2
\Big]
dy\\
& &\;\;\;\;\;\;\;=\int_{B_R \backslash \{0\}} \Big[|x|^{-\theta}
|\nabla_x (|x|^{N'-2} \psi)|^2 |x|^4-p |x|^{2(N-2)} |x|^{4+l}
|v|^{p-1} \psi^2\Big] |x|^{-2N} dx\\
& &\;\;\;\;\;\;\;=\int_{B_R \backslash \{0\}} \Big[|x|^\theta
|\nabla_x
\psi|^2-p |x|^l |v|^{p-1}\psi^2\Big] dx\\
& & \;\;\;\;\;\;\;\;\;\;\;\;\;\;\;+\int_{B_R \backslash \{0\}}
|x|^{4-\theta-2N} \Big[(N'-2)^2 |x|^{2(N'-3)} \psi^2+2(N'-2)
|x|^{2(N'-2)} \frac{x \cdot \nabla_x \psi}{|x|^2} \psi \Big] dx\\
& &\;\;\;\;\;\;=\int_{B_R \backslash \{0\}} \Big[|x|^\theta
|\nabla_x
\psi|^2-p |x|^l |v|^{p-1}\psi^2\Big] dx\\
& &\;\;\;\;\;\;\;\;\;\;\;\;\;\;\;+\int_{B_R \backslash \{0\}}
|x|^{4-\theta-2N} \Big[(N'-2)^2 |x|^{2(N'-3)} \psi^2+(N'-2)
|x|^{2(N'-3)} x \cdot \nabla_x (\psi^2) \Big] dx.
\end{eqnarray*}
Using integration by parts we find
\begin{eqnarray*}
& &\int_{B_R \backslash \{0\}} |x|^{4-\theta-2N} \Big[(N'-2)^2
|x|^{2(N'-3)} \psi^2+(N'-2) |x|^{2(N'-3)} x \cdot \nabla_x (\psi^2)
\Big]
dx\\
& &\;\;\;\;\;\;\;=(N'-2)^2 \int_{B_R \backslash \{0\}}
|x|^{\theta-2} \psi^2-(N'-2) \int_{B_R \backslash \{0\}} \mbox{div}
\Big(|x|^\theta \frac{x}{|x|^2} \Big) \psi^2 dx\\
& &\;\;\;\;\;\;\;=0.
\end{eqnarray*}
Hence
$$\int_{\R^N \backslash {\overline {B_{1/R}}}} \Big[|y|^\theta
|\nabla_y {\tilde \psi}|^2-p |y|^\beta |w|^{p-1} {\tilde \psi}^2
\Big] dy=\int_{B_R \backslash \{0\}} \Big[|x|^\theta |\nabla_x
\psi|^2-p |x|^l |v|^{p-1}\psi^2\Big] dx.$$ The proposition clearly
follows from this identity. \qed

\subsection{Proof of Proposition \ref{p5.2}}
Firstly we recall that  $\theta=-2 \sigma$, $l=\alpha-\sigma (p+1)$
and $\sigma=\frac{1}{2} [N-2-\sqrt{(N-2)^2-4 \ell}]$. Moreover
$u=|x|^{-\sigma} v \in H^1_{loc} (\Omega)$,  $v \in
H^{1,\theta}_{loc} (\Omega)\cap L_{loc}^\infty(\Omega)$.

For $\phi \in H^1_{c}(\Omega)$ with $|x|^\sigma \phi\in
L_{loc}^\infty(\Omega)$, we define ${\tilde \phi}(x)=|x|^\sigma
\phi(x)$. Then by Remark \ref{space}, ${\tilde \phi} \in
H^{1,\theta}_{c} (\Omega)\cap L_{loc}^\infty(\Omega)$. Moreover,
with $l=\alpha-\sigma (p+1)$ and $\sigma^2-(N-2) \sigma+\ell=0$, we
have
\begin{eqnarray*}
Q_v ({\tilde \phi})&=& \int_\Omega \Big[|x|^\theta |\nabla
{\tilde \phi}|^2-p|x|^l |v|^{p-1} {\tilde \phi}^2\Big]\\
&=&\int_\Omega \Big[|x|^{-2 \sigma} |\nabla {\tilde
\phi}|^2-(\sigma^2-(N-2) \sigma+\ell)|x|^{-2\sigma-2} {\tilde
\phi}^2-p |x|^\alpha |u|^{p-1} \phi^2\Big]\\
&=&\int_\Omega \Big[|x|^{-2 \sigma} |\nabla {\tilde \phi}|^2+\sigma
{\tilde \phi}^2 \nabla (|x|^{-2 \sigma-2} x)+(\sigma^2-\ell) |x|^{-2
\sigma-2} {\tilde \phi}^2-p |x|^\alpha |u|^{p-1} \phi^2\Big]\\
&=& \int_\Omega \Big[|x|^{-2 \sigma} |\nabla {\tilde \phi}|^2-2
\sigma |x|^{-2 \sigma} {\tilde \phi} \nabla {\tilde \phi} \cdot
\frac{x}{|x|^2}+(\sigma^2-\ell)|x|^{-2 \sigma-2} {\tilde \phi}^2-p |x|^\alpha |u|^{p-1} \phi^2\Big]\\
&=&\int_\Omega \Big[|\nabla (|x|^{-\sigma} {\tilde \phi})|^2-\ell
|x|^{-2} |x|^{-2 \sigma} {\tilde \phi}^2-p |x|^\alpha |u|^{p-1}
\phi^2\Big]\\
&=& \int_\Omega \Big[|\nabla \phi|^2-\ell |x|^{-2} \phi^2-p
|x|^\alpha
|u|^{p-1} \phi^2\Big] \\
&=&\mathcal{Q}_u(\phi).
\end{eqnarray*}

The conclusion of the proposition follows easily from the above identity. \qed

\subsection{Proof of Proposition \ref{L}} This follows from a simple
calculation.
 For any given $\psi \in H^1_c (B_R \backslash \{0\})\cap L_{loc}^\infty(B_R \backslash \{0\})$, we define
$${\tilde \psi} (y)= \psi (x), \;\;\;\;y=\frac{x}{|x|^2}.$$
Clearly ${\tilde \psi} \in H^1_c (\R^N \backslash {\overline
{B_{1/R}}})\cap L_{loc}^\infty(\R^N \backslash {\overline
{B_{1/R}}})$, and
\begin{eqnarray*}
& &\int_{\R^N \backslash {\overline {B_{1/R}}}} \Big[|y|^{\tilde
\theta }|\nabla_y {\tilde \psi}|^2-p |y|^{\tilde l} |z|^{p-1}
{\tilde \psi}^2 \Big]
dy\\
& &\;\;\;\;\;\;\;=\int_{B_R \backslash \{0\}}
\Big[|x|^{-\tilde\theta} |\nabla_x \psi|^2 |x|^4-p |x|^{-\tilde l}
|v|^{p-1} \psi^2\Big] |x|^{-2N} dx\\
& &\;\;\;\;\;\;\;=\int_{B_R \backslash \{0\}} \Big[|x|^\theta
|\nabla_x \psi|^2-p |x|^l |v|^{p-1}\psi^2\Big] dx.
\end{eqnarray*}
 The
conclusion of the proposition is a direct consequence of the above
identity. \qed

\section{Proof of Theorem \ref{F'}}
\setcounter{equation}{0}

We need the following lemma.

\begin{lem}
\label{l2.1}

Suppose that $p>\frac{N'+2+2 \tau}{N'-2}$, $N'>2$ and $\tau>-2$.
Then for every $\kappa>0$, problem \eqref{p*} with $\Omega=\R^N$ has a
unique positive radial solution $v_\kappa$ satisfying $v(0)=\kappa$.
Moreover, $v_\kappa$ is of the form
$$v_\kappa (r)=\kappa v_1 (\kappa^{\frac{p-1}{\tau+2}} r)$$
where $v_1$ is the unique solution of the problem
\begin{equation}
\label{2.1} \left \{ \begin{array}{l} (r^{N-1+\theta}
v'(r))'+r^{N-1+l} v^p (r)=0, \;\;\; r>0,\\
v(0)=1, \;\;\;\;\;\; \lim_{r \to 0^+} r^{N-1+\theta} v'(r)=0,
\end{array} \right.
\end{equation}
and $v_\kappa$ has the properties:

(i) for every $\kappa>0$,
\begin{equation}
\label{2.2} \lim_{r \to \infty} r^{\frac{2+\tau}{p-1}} v_\kappa
(r)=\Big \{ \frac{2+\tau}{p-1} \Big(N'-2-\frac{2+\tau}{p-1} \Big)
\Big \}^{1/(p-1)}.
\end{equation}

(ii) for $p \geq p_c (N', \tau)$,
\begin{equation}
\label{2.3} v_\kappa (r)<V_\infty (r):=C_0 r^{-\frac{2+\tau}{p-1}}
\;\;\;\; \forall r>0, \;\; \forall \kappa>0.
\end{equation}
\end{lem}

{\bf Proof.} If $\theta=0$, this is Lemma 4.1 in \cite{DDG}, which
follows from results in \cite{KYY, Li, Wa}. Since the ODE satisfied
by $u_\kappa(r)$ here is exactly the same as that satisfied by the
radial solution in Lemma 4.1 of \cite{DDG} once $(N, \alpha)$ there
is replaced by $(N', \tau)$, the conclusions here follow from the
same reasoning as in \cite{DDG} if we replace $(N, \alpha)$ there by
$(N',\tau)$.

The conclusions of this lemma are also contained in Corollary 1.3 of
\cite{DLF}, where radial solutions of more general equations are
considered. \qed

{\bf Proof of Theorem \ref{F'}.} We first show the nonexistence of
nontrivial stable solutions of \eqref{p*} for $1<p<p_c (N', \tau)$.
Arguing indirectly we assume that $1<p<p_c (N', \tau)$ and \eqref{p*}
has a solution $v \not \equiv 0$ that is stable. We are going to
deduce a contradiction.

For every $R>0$, we define the test function $\psi_R (x)=\varphi
(\frac{|x|}{R})$, where $\varphi \in C^2 (\R)$, $0 \leq \varphi \leq
1$ everywhere on $\R$ and
$$\varphi (t)=\left \{ \begin{array}{ll} 1 \;\;\;\; \mbox{if $|t|
\leq 1$}, \\
0 \;\;\;\; \mbox{if $|t| \geq 2$}.
\end{array} \right.$$
We observe that for any $\gamma \in [1, 2p+2 \sqrt{p(p-1)}-1)$ and
any $m \geq \max \{\frac{p+\gamma}{p-1}, 2\}$, Proposition \ref{K}
gives
\begin{eqnarray*}
\int_{B_R} |x|^l |v|^{p+\gamma} &\leq& C \int_{B_{2R} \backslash
B_R} |x|^{\frac{\theta (\gamma+p)-l (\gamma+1)}{p-1}} \Big[|\nabla
\psi|^2+|\psi||\Delta \psi|+|\psi| \frac{|\nabla \psi|}{|x|}
 \Big]^{\frac{\gamma+p}{p-1}}\\
&\leq& C R^{N'-\frac{(2+\tau) \gamma+2p+\tau}{p-1}} \;\;\;\;\;\;\;\;
\forall R>0.
\end{eqnarray*}
where $C$ is a positive constant independent of $R$.

Consider the function
$$\Delta (N', p, \gamma, \tau)=N' (p-1)-(2+\tau) \gamma-2p-\tau,$$
and define
$$\gamma (p)=2p+2 \sqrt{p(p-1)}-1, \;\;\; \Gamma (p)=\frac{(2+\tau)
\gamma (p)+2p+\tau}{p-1}.$$ As in the proof of Theorem 2.1 in
\cite{DDG} we can rewrite $\Gamma(p)$ in the form
\[
\Gamma(p)=2(2+\tau)\left(1+\frac{1}{p-1}+\sqrt{1+\frac{1}{p-1}}\right)+2
\]
which shows that $\Gamma(p)$ is strictly decreasing in $p$ for
$p>1$, with $\Gamma(1)=+\infty$ and $\Gamma(+\infty)=10+4\tau$.
Therefore $\Delta(N',p,\gamma(p), \tau)=(p-1)(N'-\Gamma(p))<0$ for
all $p>1$ when $N'\leq 10+4\tau$, and for $N'>10+4\tau$, there is a
unique $p^*=p^*(\tau)\in (1,\infty)$ such that $N'=\Gamma(p^*)$ and
\[
(p-p^*)\Delta(N', p,\gamma(p), \tau)>0 \mbox{ for } p\in
(1,\infty),\; p\not= p^*.
\]
We note that
 $N'=\Gamma(p^*)$ is equivalent to
 \[
 \left(\frac{N'-2}{2+\tau}-2\right)p^*-\frac{N'-2}{2+\tau}=2\sqrt{p^*(p^*-1)}.
\]
It follows that
\[
p^*>\frac{N'-2}{2+\tau}\left(\frac{N'-2}{2+\tau}-2\right)^{-1}>\frac{N'+2+2\tau}{N'-2}>\tilde p_c(N', \tau).
\]
On the other hand, a simple calculation shows that the equation
\[
 \left[\left(\frac{N'-2}{2+\tau}-2\right)p^*-\frac{N'-2}{2+\tau}\right]^2=4{p^*(p^*-1)}
\] is equivalent to
\[
a(p^*)^2-bp^*+c=0 \mbox{ with } a,\; b,\; c \mbox{ given by } \eqref{abc}.
\]
Thus we necessarily have
 $p^*=p_c(N',\tau)$, and
\begin{eqnarray*}
\Delta (N',p, \gamma (p), \tau)=0 && \mbox{for $p=p_c (N',
\tau)$};\\
\Delta (N',p, \gamma (p), \tau)<0 && \mbox{for $1<p<p_c (N',\tau)$}.
\end{eqnarray*}
Since we have assumed $1<p<p_c (N', \tau)$, we can choose $\gamma
\in (1, \gamma (p))$ close enough to $\gamma (p)$ such that
$$N'-\frac{(2+\tau) \gamma+2p+\tau}{p-1}<0.$$
Fix such a $\gamma$ and let $R \to +\infty$ in our earlier
inequality, we conclude that
$$\int_{\R^N} |x|^l |v|^{\gamma+p}=0.$$
This implies $|v|^{\gamma+p} \equiv 0$ in $\R^N$; a contradiction.

Next we show that if $p \geq p_c (N', \tau)$ (which is possible only
if $N'>10+4 \tau$), then for every $\kappa>0$, the positive radial
solution $v_\kappa$ defined in Lemma \ref{l2.1} is a stable solution
of \eqref{p*}.

We first show $v_\kappa \in H^{1, \theta}_{loc} (\R^N)$. We only
need to show that for any $R>1$,
$$
\int_{B_R} |x|^{\theta} v_\kappa^2<\infty, \;\;\;\;\;
 \int_{B_R}
|x|^\theta |\nabla v_\kappa|^2<\infty.
$$
Since $v_\kappa\in L^\infty_{loc}(\R^N)$, the first inequality is an
easy consequence of the assumption that $N'=N+\theta>2$. We now show
that
 $\int_{B_R} |x|^{\theta} |\nabla v_\kappa|^2 dx<\infty$. It follows from the equation of $v_\kappa$ that $v_\kappa'(r)<0$
for $r>0$. Moreover,
\begin{eqnarray*}
|v_\kappa'(r)| &= &r^{1-N-\theta} \int_0^r s^{N-1+l} v_\kappa^p (s)
ds\\
&\leq& r^{1-N-\theta} \int_0^r s^{N-1+l} V_\infty^p (s)
ds \;\;\;\;\; \mbox{by \eqref{2.3}}\\
&=& C_0^p r^{1-N-\theta} \int_0^r s^{N-1+l-\frac{p (2+\tau)}{p-1}}
ds\\
&=& C_{p, N', \tau} r^{1+\tau-\frac{p (2+\tau)}{p-1}} \;\;\;
(\mbox{note that $N-1+l-\frac{p (2+\tau)}{p-1}>-1$ for
$p>\frac{N'+\tau}{N'-2}$}).
\end{eqnarray*}
Therefore, for any $R>0$ and $p \geq p_c (N', \tau)\; (>\frac{N'+2+2
\tau}{N'-2})$, we have $N+1+\theta+2 \tau-\frac{2p
(2+\tau)}{p-1}>-1$ and
$$\int_{B_R} |x|^\theta |\nabla v_\kappa|^2=\int_0^R r^{N-1+\theta} |v'(r)|^2 dr
 \leq C_{p, N', \tau}^2  \int_0^R r^{N+1+\theta+2 \tau-\frac{2p (2+\tau)}{p-1}}
 dr<\infty.$$

Since \eqref{2.3} holds, we have, for every $\psi \in C_0^1 (\R^N)$,
\begin{eqnarray*}
Q_{v_\kappa} (\psi) &=& \int_{\R^N} |x|^\theta |\nabla \psi|^2-p
\int_{\R^N} |x|^l v_\kappa^{p-1} \psi^2\\
&\geq& \int_{\R^N} |x|^\theta |\nabla \psi|^2-p \int_{\R^N} |x|^l
V_\infty^{p-1} \psi^2 \\
&=& \int_{\R^N} |x|^\theta |\nabla \psi|^2- \int_{\R^N} p C_0^{p-1}
|x|^l |x|^{-(2+\tau)} \psi^2\\
&=&  \int_{\R^N} |x|^\theta |\nabla \psi|^2- \int_{\R^N} p C_0^{p-1}
|x|^{-(2-\theta)} \psi^2.
\end{eqnarray*}
By the  Caffarelli-Kohn-Nirenberg inequality \cite{CKN},
$$ \Big(\int_{\R^N} \frac{|\psi|^q}{|x|^{b q}} dx \Big)^{2/q} \leq C
(N,a,b) \int_{\R^N} \frac{|\nabla \psi|^2}{|x|^{2a}} dx,$$ where
$C(N,a,b)$ is a positive constant and
$$
-\infty<a<\frac{N-2}{2}, \;\;\; a \leq b \leq a+1, \;\;\;
q=\frac{2N}{N-2+2(b-a)}.
$$
In our case here,
$$a=-\frac{\theta}{2}, \;\;\;\; b=1-\frac{\theta}{2}=1+a, \;\;\; \; q=2,$$
and by \cite{CW}, $C(N,-\frac{\theta}{2},1-\frac{\theta}{2})$ has
the optimal value $\frac{4}{(N'-2)^2}$. Therefore
$$\int_{\R^N} \frac{|\psi|^2}{|x|^{2-\theta}} dx \leq
\frac{4}{(N'-2)^2} \int_{\R^N} |x|^\theta |\nabla \psi|^2 dx,
$$
and
\begin{equation}
\label{add} \int_{\R^N} |x|^\theta |\nabla \psi|^2- \int_{\R^N} p
C_0^{p-1} |x|^{-(2-\theta)} \psi^2 \geq \Big(\frac{(N'-2)^2}{4}-p
C_0^{p-1} \Big) \int_{\R^N} |x|^{-(2-\theta)} \psi^2 \geq 0,
\end{equation}
since
$$
\frac{(N'-2)^2}{4}-p C_0^{p-1}=\frac{(N'-2)^2}{4}-f(p) \geq 0 \;\;\;
\mbox{for $p \geq p_c (N', \tau)$}.
$$
Thus $Q_{v_\kappa} (\psi) \geq 0$. This means that
$v_\kappa$ is a stable solution of \eqref{p*}. This completes the
proof. \qed

\section{Asymptotic bounds  and related results}
\setcounter{equation}{0}

In this section, we supply the proofs of Theorems \ref{G} and
\ref{H}, and also prove the necessity of the assumption $\tau>-2$
and the uniqueness of the radial solution $V_\infty$.
\subsection{Proof of Theorem \ref{G}}
 Since $v$ has finite Morse index, it is stable outside
a compact subset of $\Omega$ and hence there exists $R_*>0$ small
such that $v$ is stable in $B_{R_*} \backslash \{0\}$.

{\it Step 1.} Suppose that $v$ is a stable positive solution of
\eqref{p*} in $B_{R_*} \backslash \{0\}$. Then for every $\gamma \in
[1, 2p+2 \sqrt{p(p-1)}-1)$ and every open ball $B_R (y)$ with
$0<|y|<\frac{4}{5} R_*$ and $R=\frac{|y|}{4}$, we have
\begin{equation}
\label{3.1} \int_{B_R (y)} |x|^l v^{\gamma+p} \leq C
R^{N'-\frac{(2+\tau) \gamma+2p+\tau}{p-1}},
\end{equation}
where $C$ is a positive constant depending on $m, p, N', \tau$ but
not on $y$.

Since $v$ is stable in $B_{R_*} \backslash \{0\}$, Proposition
\ref{K} holds when $\Omega=B_{R_*} \backslash \{0\}$. We fix a
function $\varphi_0 \in C^2 (\R)$ satisfying $0 \leq \varphi_0 \leq
1$ everywhere on $\R$ and
$$\varphi_0 (t)=\left \{ \begin{array}{ll} 0 \;\;\;\; \mbox{if $|t|
\leq 1$}, \\
1 \;\;\;\; \mbox{if $|t| \geq 2$}.
\end{array} \right.$$
We then apply Proposition \ref{K} with $m=1+\max
\{\frac{p+\gamma}{p-1},2\}$ and test function $\psi (x):=\varphi_0
(\frac{|x-y|}{R})$ and obtain \eqref{3.1} as we did in the proof of
Theorem \ref{F'}.

{\it Step 2.} Suppose that $v$ is a stable solution of \eqref{p*} in
$B_{R_*} \backslash \{0\}$. Then if $1<p<p_c (N,0)$, there exists a
small $\epsilon_0=\epsilon_0 (p)>0$ such that for every $\epsilon
\in [0, \epsilon_0]$ and every open ball $B_{2R} (y)$ with $0<|y|
\leq \frac{2}{3} R_*$ and $R=|y|/8$, we have
\begin{equation}
\label{3.2} R^{-\frac{N \theta}{2-\epsilon}} \int_{B_{2R} (y)}
\Big(|x|^l v^{p-1} \Big)^{\frac{N}{2-\epsilon}} \leq C
R^{N-\frac{2N}{2-\epsilon}},
\end{equation}
where $C$ is a positive constant depending on $m,p,N, \tau$ but not
on $y$ and $\epsilon$.

Let us recall that for
$$\Delta (N', p, \gamma, \tau)=N' (p-1)-(2+\tau) \gamma-2p-\tau$$
and
$$\gamma (p)=2p+2 \sqrt{p(p-1)}-1, $$ we have
\begin{eqnarray*}
\Delta (N',p, \gamma (p), \tau)=0 && \mbox{for $p=p_c (N',
\tau)$};\\
\Delta (N',p, \gamma (p), \tau)<0 && \mbox{for $1<p<p_c (N',\tau)$}.
\end{eqnarray*}
Taking $\tau=0$ we obtain
\begin{equation}
\label{3.3} \Delta (N,p, \gamma (p),0)=N (p-1)-2(\gamma (p)+p)<0
\;\;\;\; \mbox{for $1<p<p_c (N,0)$}.
\end{equation}
Thus we can fix $\gamma_*=\gamma_* (p) \in (1, \gamma (p))$ such
that
\begin{equation}
\label{3.4} \frac{p+\gamma_*}{(p-1) N/2}>1.
\end{equation}
It is seen from \eqref{3.4} that we can find $\epsilon_0=\epsilon_0
(p)>0$ sufficiently small so that
$$\frac{p+\gamma_*}{(p-1) \rho}>1 \;\;\;\; \forall \rho \in
\Big[\frac{N}{2}, \frac{N}{2-\epsilon_0} \Big].$$ Fix such a $\rho$
and set
$$\xi=\frac{p+\gamma_*}{(p-1) \rho}.$$
By H\"older's inequality and \eqref{3.1},
\begin{eqnarray*}
\int_{B_{2R} (y)} (|x|^l v^{p-1})^\rho &\leq& \Big(\int_{B_{2R} (y)}
|x|^l v^{\gamma_*+p} \Big)^{1/\xi} \Big(\int_{B_{2R} (y)}
|x|^{\frac{l (\rho \xi-1)}{\xi-1}} \Big)^{(\xi-1)/\xi}\\
&\leq& C R^{(N'-\frac{(2+\tau) \gamma_*+2p+\tau}{p-1})
\frac{1}{\xi}}
R^{(N+\frac{l (\rho \xi-1)}{\xi-1}) \frac{\xi-1}{\xi}}\\
&=&C R^{N-2 \rho+\frac{(p+\gamma_*) \theta}{(p-1) \xi}},
\end{eqnarray*}
which implies that
\begin{equation}
\label{3.6} R^{-\theta \rho} \int_{B_{2R} (y)} (|x|^l v^{p-1})^\rho
\leq C R^{N-2 \rho},
\end{equation}
and  \eqref{3.2} follows if we take $\rho=\frac{N}{2-\epsilon}$.

{\it Step 3.} Harnack inequality: Under the conditions of Step 2,
there exists a positive constant $K$ such that
\begin{equation}
\label{3.10} \max_{|x|=r} v(x) \leq K \min_{|x|=r} v(x) \;\;\;\;
\forall r \in (0, R_*].
\end{equation}

Regarding $v=v(x)$ as a solution of the equation
$$\mbox{div} (|x|^\theta \nabla v)+d(x)v=0$$
with $d(x)=|x|^l v^{p-1} (x)$, in view of \eqref{3.2}, we can apply
Harnack's inequality on each ball $B_R (y)$ with $0<|y|<\frac{2}{3}
R_*$, $R=\frac{|y|}{8}$, to obtain
\begin{equation}
\label{3.11} \sup_{B_R (y)} v \leq K \inf_{B_R (y)} v,
\end{equation}
where $K$ depends on $N',m, p, \tau$ and $R^\epsilon \|R^{-\theta}
d\|_{L^{\frac{N}{2-\epsilon}} (B_{2R} (y))}$ (see \cite{GT} p. 209).
(Note that for $x \in B_{2R} (y)$, $|y|-|x-y| \leq |x| \leq
|y|+|x-y|$ and thus $6 R \leq |x| \leq 10 R$. This implies that
$|x|^\theta \geq 6^\theta R^\theta$ provided $\theta \geq 0$;
$|x|^\theta \geq 10^\theta R^\theta$ provided $\theta<0$. Therefore,
the $\lambda$ in \cite{GT} is $6^\theta R^\theta$ or $10^\theta
R^\theta$.) Due to \eqref{3.2},
$$R^\epsilon \|R^{-\theta} d\|_{L^{\frac{N}{2-\epsilon}} (B_{2R}
(y))} \leq R^\epsilon C R^{-\epsilon}=C.$$ Therefore, $K$ is
independent of $R$. Given any $r \in (0, \frac{2}{3} R_*]$, the
sphere $\{|x|=r\}$ can be covered by a finite number of balls of the
form $B_R (y)$ with $|y|=r$ and $R=|y|/8=r/8$, and this finite
number is independent of $r$. Therefore, by enlarging $K$ in
\eqref{3.11} properly, we have
$$\max_{|x|=r} v(x) \leq K \min_{|x|=r} v(x) \;\;\;\;
\forall r \in \Big(0, \frac{2}{3} R_* \Big].$$ Since $v$ is positive
and continuous in $\{\frac{2}{3} R_* \leq |x| \leq R_*\}$, by
further enlarging $K$ if necessary, we can guarantee that the above
inequality holds for all $r \in (0, R_*]$, and \eqref{3.10} is
proved.

{\it Step 4.} Under the conditions of Step 2, there exists a
positive constant $C$ such that
\begin{equation}
\label{3.12} v(x) \leq C |x|^{-\frac{2+\tau}{p-1}} \;\;\;\; \forall
x \in B_{R_*} \backslash \{0\}.
\end{equation}
From \eqref{3.2} with $\epsilon=0$ we obtain, for
$0<|y|<\frac{2}{3}R_*$ and $R=|y|/8$,
$$\vartheta \Big[\inf_{B_R (y)} v \Big]^{\frac{N(p-1)}{2}}
R^{N(1+\frac{\tau}{2})} \leq \int_{B_R (y)} \Big[R^{-\theta} |x|^l
v^{p-1} \Big]^{\frac{N}{2}} \leq C,$$ where $\vartheta:=\vartheta
(N)$ is a positive constant independent of $y$ and $v$. It follows
that
$$\inf_{B_R (y)} v \leq \Big(\frac{C}{\vartheta}
\Big)^{\frac{2}{N(p-1)}} R^{-\frac{2+\tau}{p-1}}.$$ We can now apply
\eqref{3.10} to obtain
$$\sup_{B_R (y)} v \leq K \Big(\frac{C}{\vartheta}
\Big)^{\frac{2}{N(p-1)}} R^{-\frac{2+\tau}{p-1}}.$$ In particular,
$$v (y) \leq K \Big(\frac{C}{\vartheta}
\Big)^{\frac{2}{N(p-1)}} R^{-\frac{2+\tau}{p-1}}=C_1
|y|^{-\frac{2+\tau}{p-1}}$$ for all $y$ satisfying $0<|y| \leq
\frac{2}{3} R_*$. Since both $v(y)$ and $|y|^{-\frac{2+\tau}{p-1}}$
are positive and continuous on $\{\frac{2}{3} R_* \leq |y| \leq
R_*\}$, by enlarging $C_1$ if necessary, we have
\begin{equation}
\label{3.13} v(y) \leq C_1 |y|^{-\frac{2+\tau}{p-1}} \;\;\;
\mbox{for all $y$ satisfying $0<|y| \leq R_*$}.
\end{equation}
The proof is complete. \qed

\begin{rem}
\label{r3.1} The condition $1<p<p_c (N,0)$ in Theorem \ref{G} is
only used to obtain \eqref{3.4}. We may attempt to  replace it by
other conditions. For example, since $\Delta (N',p, \gamma (p),
\tau)<0$ for $1<p<p_c (N',\tau)$, we see that
\begin{equation}
\label{3.3-1} N (p-1)-2 (\gamma (p)+p)<(\gamma (p)+1)
\Big[\tau-\frac{(p-1)}{(1+\gamma (p))} \theta\Big] \leq 0 \;\;\;
\mbox{for $1<p<p_c (N',\tau)$}
\end{equation}
provided $\tau-\frac{(p-1)}{(1+\gamma (p))} \theta \leq 0$.
Therefore, we can fix $\gamma_* \in [1, \gamma (p))$ such that
\eqref{3.4} holds provided
$$\tau \leq \frac{(p-1)}{2p+2
\sqrt{p(p-1)}} \theta \;\;\; \mbox{and} \;\;\; 1<p<p_c (N', \tau).
\leqno(A)$$ However, it is easy to see that condition (A) is more
restrictive than requiring $1<p<p_c(N', 0)$, because we will show
below that
 the
function $p_c (N', \cdot)$ is increasing for any fixed $N'$, and
thus $\tau \leq \frac{(p-1)}{2p+2 \sqrt{p(p-1)}} \theta$ implies
\begin{equation}
\label{3.3-3} p_c (N', \tau) \leq p_c \Big(N', \frac{(p-1)}{2p+2
\sqrt{p(p-1)}} \theta \Big)=p_c (N,0).
\end{equation}
To see the equality above, we note that if ${\tilde
\tau}=\frac{(p-1)}{2p+2 \sqrt{p(p-1)}} \theta$, then
\begin{eqnarray*}
& & \Delta (N', p, \gamma (p), {\tilde \tau})\\
& & \;\;\;\;\;\;\;\;\;\;\;=N (p-1)-2 (p+\gamma (p))-(1+\gamma (p))
\Big({\tilde \tau}- \frac{(p-1)}{1+\gamma (p)} \theta \Big )\\
& &\;\;\;\;\;\;\;\;\;\;\;= N (p-1)-2 (p+\gamma (p))\\
& &\;\;\;\;\;\;\;\;\;\;\;=\Delta (N', p, \gamma(p), 0).
\end{eqnarray*}
Hence from $\Delta (N', p_c (N', {\tilde \tau}), \gamma (p_c (N',
{\tilde \tau})), {\tilde \tau})=0$ we obtain
$$N [ p_c(N', {\tilde \tau})-1]-2 [p_c (N', {\tilde \tau})+\gamma
(p_c (N', {\tilde \tau}))]=0$$ and thus $p_c (N', {\tilde \tau})=p_c
(N,0)$.

On the other hand, if
$$\tau>\tilde\tau=\frac{(p-1)}{2p+2
\sqrt{p(p-1)}} \theta \;\;\; \mbox{and} \;\;\; 1<p<p_c (N', \tau),
\leqno(B)$$ then
$$p_c (N', \tau)>p_c (N', \tilde \tau) =p_c (N,0).$$

We now show that $p_c (N', \tau)$ is decreasing in $N'$ and
increasing in $\tau$. Recall that, for $N'>4\tau+10$,
$p_c(N',\tau)\in (1,\infty)$ is the unique solution of
\[
N'=\Gamma(p)=2(2+\tau)\left(1+\frac{1}{p-1}+\sqrt{1+\frac{1}{p-1}}\;\right)+2,
\]
which is equivalent to
\[
\frac{N'-2}{2+\tau}=2\left(1+\frac{1}{p-1}+\sqrt{1+\frac{1}{p-1}}\;\right).
\]
Since the term on the left hand side is increasing in $N'$ and
decreasing in $\tau$, while the term on the right hand side is a
decreasing function of $p$, it follows immediately that
$p_c(N',\tau)$ is increasing in $\tau$ and decreasing in $N'$.
\end{rem}

\subsection{Proof of Theorem \ref{H}}
 Since $v$ has finite Morse index, it is stable outside
a compact subset of $\Omega$ and hence there exists $R_*>0$ large
such that $v$ is stable in $\R^N \backslash {\overline {B_{R_*}}}$.

Define
$$w(y)=|x|^{N'-2} v(x), \;\;\;\;\; y=\frac{x}{|x|^2}.$$
Then $w$ satisfies
\begin{equation}
\label{4.2} -\mbox{div}(|y|^\theta \nabla w)=|y|^\beta w^p
\;\;\;\mbox{in $B_{1/R_*} \backslash \{0\}$},
\end{equation}
with
\[\mbox{
$\tau':=\beta-\theta =(N'-2)(p-1)-(4+\tau)>-2$ if
$p>(N'+\tau)/(N'-2)$.}
\]

 By Proposition \ref{F}, $w$ is a stable
positive solution of \eqref{4.2}. Therefore when $p\in \left(\frac{N'+\tau}{N'-2}, p_c(N',0)\right)$,
we can apply Theorem \ref{G} to \eqref{4.2} to conclude that
$$|y|^{\frac{2+\beta-\theta}{p-1}} w(y) \leq C \;\;\;\; \mbox{for
all small $|y|>0$},$$ which is equivalent to
$$
|x|^{\frac{2+\tau}{p-1}} v(x) \leq C \;\;\;\; \mbox{for
all large $|x|>0$}.
$$

It remains to consider the case that $p\in \left(1,
\frac{N'+\tau}{N'-2}\right]$, which implies that $\tau' \leq -2$. By
Theorem \ref{t3.3} below, in this case, \eqref{4.2} does not have a
positive solution over any punctured ball $B_R\setminus\{0\}$, which
implies that \eqref{p*} has no positive solution over any exterior
domain. Therefore there is nothing to prove for this case. \qed

\subsection{Related results}
The next result reveals the role played by the condition $\tau>-2$.

\begin{thm}
\label{t3.3} For $N'\geq 2$ and $\tau \leq -2$, problem \eqref{p*} does not admit a
positive solution over any punctured ball $B_R \backslash \{0\}
\subset \R^N \; (N \geq 2)$.
\end{thm}

{\bf Proof.} We argue indirectly by assuming that $u \in C^2 (B_R
\backslash \{0\})$ is a positive solution of \eqref{p*}. Using
spherical coordinates to write $v(x)=v(r, \omega)$ with $r=|x|$ and
$\omega=\frac{x}{|x|}$, we have
$$v_{rr}+\frac{N'-1}{r} v_r+\frac{1}{r^2} \Delta_{S^{N-1}}
v=-r^\tau v^p.$$ This equation is exactly the same as that in Theorem
2.3 of \cite{DDG} when $(N, \alpha)$ there is replaced by $(N',
\tau)$ here. Since $N'\geq 2$, the arguments  in
the proof of Theorem 2.3 in \cite{DDG} lead to a contradiction. The
proof is thus complete. \qed

Similarly, a positive radial solution $v (r)$ of \eqref{p*}
satisfies
$$v_{rr}+\frac{N'-1}{r}v_r=-r^\tau v^p,$$
which is exactly the same as that satisfied by $u (r)$ in Theorem 2.4 of
\cite{DDG} with $(N,\alpha)$ there being replaced by $(N', \tau)$
here. Thus we have the following analogue of Theorem 2.4 of \cite{DDG}.

\begin{thm}
\label{t3.4} Let $v=v(r)$ be a positive radial solution of \eqref{p*}
over $B_R \backslash \{0\}$ with ${\overline {\lim}}_{r \to 0}
v(r)=\infty$ and $p>\frac{N'+2+2 \tau}{N'-2}$. Then
$$v(r) \equiv V_\infty (r).$$
\end{thm}

This theorem implies that $V_\infty (r)$ is the unique positive
radial singular solution of \eqref{p*} over any $B_R$ when
$p>\frac{N'+2+2 \tau}{N'-2}$.

\section{Exact asymptotic behavior}

In this section, we prove Theorems \ref{I} and \ref{J}. We first prove the results for $p>\frac{N'+2+2\tau}{N'-2}$.
Then we make use of the Kelvin transformation to cover the full range of $p$.
\begin{thm}
\label{t3.2} Let $\Omega_0$ be a bounded domain in $\R^N \; (N \geq
2)$ containing 0, and let $v$ be a positive solution of \eqref{p*} with
$\Omega=\Omega_0 \backslash \{0\}$. If $v$ has finite Morse index,
then $x=0$ must be a removable singularity of $v$ provided that
\begin{equation}
\label{half-cond-p}
\frac{N'+2+2 \tau}{N'-2}< p<\min\{p_c (N',\tau),p_c(N,0)\}.
\end{equation}

On the other hand, if $p \geq p_c (N', \tau)$, then problem \eqref{p*}
has a stable positive solution with an isolated singularity at 0.
\end{thm}

{\bf Proof.} A direct calculation shows that, as long as $N'>2$ and $p>\frac{N'+\tau}{N'-2}$,
\begin{equation}
\label{3.14} V_\infty (x):=C_0 |x|^{-\frac{2+\tau}{p-1}}, \;\;\;\;
C_0=\Big \{ \frac{2+\tau}{p-1} \Big(N'-2-\frac{2+\tau}{p-1} \Big)
\Big \}^{1/(p-1)}
\end{equation}
is a positive solution of \eqref{p*} in $\R^N \backslash \{0\}$, with 0
an isolated singularity.

Moreover, when $p\geq p_c(N',\tau)$, it follows from \eqref{add} that for every $\psi \in C_0^1
(\R^N)$,
$$Q_{V_\infty} (\psi)=\int_{\R^N} [|x|^\theta |\nabla \psi|^2-p
|x|^l V_\infty^{p-1} \psi^2] \geq 0,$$ that is, $V_\infty$ is a
stable solution of \eqref{p*} on $\R^N \backslash \{0\}$. In
particular, it is a stable positive solution of \eqref{p*} in $\Omega$.

Next we suppose that \eqref{half-cond-p} holds and that $v$ is a
positive solution of \eqref{p*} with finite Morse index.
For $p$ in this range, Theorem \ref{G} applies and hence there
exist $C>0$ and small $r_0>0$ such that
\begin{equation}
\label{3.15} |x|^{\frac{2+\tau}{p-1}} v(x) \leq C \;\;\;\; \mbox{for
$0<|x|<r_0$}.
\end{equation}
Hence we can apply Theorem \ref{C} to conclude that $v$ either has a
removable singularity at $x=0$ or
\begin{equation}
\label{3.16} C_1 \leq |x|^{\frac{2+\tau}{p-1}} v(x) \leq C_2
\end{equation}
for some $C_1, C_2>0$ and small positive $|x|$, say $0<|x|<R_0$.
Thus, to complete the proof, it suffices to show that \eqref{3.16}
does not hold.

Arguing indirectly, we suppose that \eqref{3.16} holds, and then
derive a contradiction. Since $v$ has finite Morse index, we may
assume that $v$ is stable in $B_{R_*} \backslash \{0\}$ for some
sufficiently small $R_*>0$. We divide our arguments below into two
steps.

{\it Step 1.} Suppose that $N'>2$, $\tau>-2$, $p>1$ and $v$ is a
stable positive solution of \eqref{p*} in $B_{R_*} \backslash \{0\}$.
Then there exists $R_0 \in (0, R_*)$ such that for every $\gamma \in
[1, 2p+2 \sqrt{p(p-1)}-1)$ and every $r \in (0, R_0/2)$, we have
\begin{equation}
\label{3.17} \int_{r<|x|<R_0} |x|^l v^{\gamma+p} \leq C+D
r^{N'-\frac{(2+\tau) \gamma+2p+\tau}{p-1}}
\end{equation}
where $C$ and $D$ are positive constants depending on
$m,p,N',\tau,R_0,R_*$ but not on $r$.

Since $v$ is stable in $B_{R_*} \backslash \{0\}$, Proposition
\ref{K} holds with $\Omega=B_{R_*} \backslash \{0\}$. To choose a
suitable test function for our purpose here, we fix a function
$\varphi_0 \in C^2 (\R)$ as in the proof of Theorem \ref{G} and
choose another function $\varrho_0$ such that $\varrho_0 \in
C^2(\R)$, $0 \leq \varrho_0 \leq 1$ everywhere on $\R$ and
$$\varrho_0 (t)=\left \{ \begin{array}{ll} 1 \;\;\;\;& \mbox{if $t \leq
R_0$},\\ 0 \;\;\;\;& \mbox{if $t \geq (R_0+R_*)/2$}.
\end{array} \right.$$
For every $r \in (0, R_0/2)$, we define $\xi_r$ as follows
$$\xi_r (x)=\left \{ \begin{array}{ll} \varrho_0 (|x|) \;\;\;\;
& \mbox{if $|x| \geq R_0/2$},\\
\varphi_0 (\frac{2 |x|}{r}) \;\;\;\; &\mbox{if $|x| \leq R_0/2$}.
\end{array} \right.$$
Clearly $\xi_r$ belongs to $C^2_0 (B_{R_*} \backslash \{0\})$ and
satisfies $0 \leq \xi_r \leq 1$ everywhere on $\R^N$. We now choose
$m=1+\max \{\frac{p+\gamma}{p-1},2\}$ and apply Proposition \ref{K}
with $\Omega=B_{R_*} \backslash \{0\}$ and $\psi=\xi_r$ to obtain
\begin{eqnarray*}
& & \int_{r/2<|x|<R_0} |x|^l v^{\gamma+p}\\
& & \;\;\;\;\;\leq C \int_{\R^N} |x|^{\frac{\theta (p+\gamma)-l
(\gamma+1)}{p-1}} \Big(|\nabla \xi_r|^2+|\xi_r| |\Delta
\xi_r|+|\xi_r|
\frac{|\nabla \xi_r|}{|x|} \Big)^{\frac{p+\gamma}{p-1}}\\
& &\;\;\;\;\;\leq {\hat C} \Big[ \int_{R_0 \leq |x| \leq R_*}
|x|^{\frac{\theta (p+\gamma)-l (\gamma+1)}{p-1}}
\Big(|\varrho_0'(|x|)|^2+\varrho_0 (|x|) |\varrho_0''
(|x|)|+\frac{|\varrho_0'(|x|)|}{|x|} \Big)^{\frac{p+\gamma}{p-1}}\\
& &\;\;\;\;\;\;\;\;\;\;\;\;\;\;\;\;\;\;\;\;+\int_{\frac{r}{2} \leq
|x| \leq r} |x|^{\frac{\theta
(p+\gamma)-l (\gamma+1)}{p-1}} \Big(r^{-2} |\varphi_0'(2|x|/r)|^2\\
& &\;\;\;\;\;\;\;\;\;\;\;\;\;\;\;\;\;\;\;\;\;\;\;\;\;\;\;+ r^{-2}
\varphi_0 (2|x|/r) |\varphi_0'' (2 |x|/r)|+2 r^{-2}
|\varphi_0'(2|x|/r)| \Big)^{\frac{p+\gamma}{p-1}}
\Big]\\
& & \;\;\;\;\leq C_1+C_2 r^{N'-\frac{(2+\tau) \gamma+2p+\tau}{p-1}}
\end{eqnarray*}
for all $r \in (0, R_0/2)$ and all $\gamma \in [1, 2p+2
\sqrt{p(p-1)}-1)$. Hence the desired integral estimate \eqref{3.17}
holds.

{\it Step 2.} Reaching a contradiction when \eqref{half-cond-p}
holds.

Recall that
\begin{equation}
\label{3.18} \Delta (N',p, \gamma (p), \tau)<0 \;\;\; \mbox{for
$1<p<p_c (N', \tau)$.}
\end{equation}
Hence
\begin{equation}
\label{3.18-1} \Delta (N',p, \gamma (p), \tau)<0 \;\;\; \mbox{ when \eqref{half-cond-p} holds}.
\end{equation}
On the other hand,
\begin{equation}
\label{3.19} \Delta (N',p,1,\tau)=(N'-2)p-(N'+2)-2 \tau \geq 0
\;\;\; \mbox{if $p \geq \frac{N'+2+2 \tau}{N'-2}$}.
\end{equation}
Therefore, under our assumption on $p$, we can find $\gamma_0 \in
[1, \gamma (p))$ such that $\Delta (N', p,
\gamma_0, \tau)=0$, that is,
$$N'-\frac{(2+\tau) \gamma_0+2p+\tau}{p-1}=0.$$
Choosing $\gamma=\gamma_0$ in \eqref{3.17}, we obtain
$$\int_{\{r<|x|<R_0 \}} |x|^l v^{\gamma_0+p} \leq C+D.$$
On the other hand, using \eqref{3.16} we deduce
$$\int_{\{r<|x|<R_0 \}} |x|^l v^{\gamma_0+p} \geq C_1^{\gamma_0+p}
\int_{\{r<|x|<R_0\}} |x|^{l-\frac{2+\tau}{p-1}
(\gamma_0+p)}=C_1^{p+\gamma_0} \int_r^{R_0} s^{-1} ds$$
$$=C_1^{\gamma_0+p} \log (R_0/r) \to \infty \;\;\; \mbox{as $r \to
0^+$},$$ a contradiction. This completes the proof.
\qed

\begin{thm}
\label{t4.3} Suppose that $\Omega_0$ is a bounded domain containing
0 and the condition \eqref{half-cond-p} in Theorem \ref{t3.2} holds.
If $v$ is a positive solution of \eqref{p*} in $\Omega:=\R^N \backslash
\Omega_0$ that has finite Morse index, then it must be a fast decay
solution.

On the other hand, if $p \geq p_c (N', \tau)$, then \eqref{p*} admits a
stable positive solution decaying at the slower rate
$|x|^{-\frac{2+\tau}{p-1}}$ at infinity.
\end{thm}

{\bf Proof.} If $p \geq p_c (N', \tau)$, we already know from the
proof of Theorem \ref{t3.2} that $V_\infty$ is a stable positive
solution of \eqref{p*} over $\Omega$ with slow decay at infinity.

Next we suppose that \eqref{half-cond-p} holds and $v$ is a positive
solution of \eqref{p*} with finite Morse index.
Therefore, Theorem \ref{H} applies and there exists $C>0$
and large $R_*>0$ such that
\begin{equation}
\label{4.3} |x|^{\frac{2+\tau}{p-1}} v(x) \leq C \;\;\;\; \mbox{for
$|x|>R_*$}.
\end{equation}
Hence we can apply Theorem \ref{D} to conclude that either $v$ has
fast decay at infinity, or there exist $C_1$, $C_2>0$ such that
\begin{equation}
\label{4.4} C_1 \leq |x|^{\frac{2+\tau}{p-1}} v(x) \leq C_2
\;\;\;\mbox{for all large $|x|$}.
\end{equation}

Thus to complete the proof, we only have to show that \eqref{4.4}
does not hold. Suppose that \eqref{4.4} holds, we will derive a
contradiction.

Since $v$ has finite Morse index over $\Omega$, we may assume that
$v$ is stable in $\R^N \backslash B_R$.

{\it Step 1.} Suppose that $\tau>-2$, $p>1$ and $v$ is a stable
positive solution of \eqref{p*} in $\R^N \backslash B_R$ with $R>R_*$.
Then there exists $R_0>R$ such that for every $\gamma \in [1, 2p+2
\sqrt{p(p-1)}-1)$ and every $r>R_0$, we have
\begin{equation}
\label{4.5} \int_{\{R_0<|x|<r\}} |x|^l v^{\gamma+p} \leq C+D
r^{N'-\frac{(2+\tau) \gamma+2p+\tau}{p-1}},
\end{equation}
where $C$ and $D$ are positive constants depending on $m,p,N',\tau,
R, R_0$ but not on $r$.

Since $v$ is stable in $\R^N \backslash B_R$, Proposition \ref{K}
holds with $\Omega=\R^N \backslash B_R$. We now choose a suitable
test function. We fix $\varphi_0 \in C^2 (\R)$ as in the proof of
Theorem \ref{t3.2}. Then define
$${\tilde \xi}_r (x)=\left \{ \begin{array}{ll} 1
\;\;\;& \mbox{if $|x| \leq R_*/2$},\\
1-\varphi_0 (\frac{2 |x|}{r}) \;\;\; &\mbox{if $|x| \geq R_*/2$}.
\end{array} \right.$$
We may then prove \eqref{4.5} in the same way as in Step 1 of the
proof of Theorem \ref{t3.2}.

{\it Step 2.} Reaching a contradiction when \eqref{half-cond-p}
holds.

 As in the proof of Theorem \ref{t3.2}, under our assumption on
$p$, we can find $\gamma_0 \in [1, \gamma (p))$ such that $\Delta
(N', p, \gamma_0, \tau)=0$, that is
$$N'-\frac{(2+\tau) \gamma_0+2p+\tau}{p-1}=0.$$
Choosing $\gamma=\gamma_0$ in \eqref{4.5}, we obtain
$$\int_{\{R_0<|x|<r\}} |x|^l v^{\gamma_0+p} \leq C+D.$$
On the other hand, using \eqref{4.4} we deduce
$$\int_{\{R_0<|x|<r \}} |x|^l v^{\gamma_0+p} \geq C_1^{\gamma_0+p}
\int_{\{R_0<|x|<r\}} |x|^{l-\frac{2+\tau}{p-1}
(\gamma_0+p)}=C_1^{p+\gamma_0} \int_{R_0}^r s^{-1} ds$$
$$=C_1^{\gamma_0+p} \log (r/R_0) \to \infty \;\;\; \mbox{as $r \to
\infty$},$$ a contradiction. This completes our proof.
\qed

We next use the Kelvin transformation to show that the conclusions of both Theorems \ref{t3.2} and
\ref{t4.3} continue to hold when
\begin{equation}
\label{4.6} \tilde p_c(N', \tau)<p<\min\left\{\frac{N'+2+2 \tau}{N'-2},p_c(N,0)\right\};
\end{equation}
and moreover, when $p\in \left(\frac{N'+\tau}{N'-2}, \tilde
p_c(N',\tau)\right]$, $V_\infty$ is a stable solution of \eqref{p*}
over $\R^N\setminus\{0\}$. Clearly Theorems \ref{I} and \ref{J}
follow from these.

We only consider the case of Theorem \ref{t3.2}, the proof for the
case of Theorem \ref{t4.3} is analogous.

\begin{thm}
\label{t4.4} Let $\Omega_0$ be a bounded domain in $\R^N \; (N \geq
2)$ containing 0, and let $v$ be a positive solution of \eqref{p*} with
$\Omega=\Omega_0 \backslash \{0\}$. If $v$ has finite Morse index,
then $x=0$ must be a removable singularity of $v$ provided that \eqref{4.6} holds.

On the other hand, if  $\frac{N'+\tau}{N'-2}< p \leq \tilde p_c(N',
\tau)$, then problem \eqref{p*} has a stable positive solution with
an isolated singularity at 0.
\end{thm}

{\bf Proof.}\
Let $p$ be in the range given by \eqref{4.6} and suppose that $v$ is a positive
solution of \eqref{p*} with finite Morse index. Then there exists $R>0$
such that $v$ is stable in $B_R \backslash \{0\}$. Therefore, the
function $w$ given by
$$w(y)=|x|^{N'-2} v(x), \;\;\; y=\frac{x}{|x|^2}$$
is a stable solution of
\begin{equation}
\label{4.7} -\mbox{div} (|y|^\theta \nabla w)=|y|^\beta w^p \;\;\;\;
\mbox{in $\R^N \backslash {\overline {B_{1/R}}}$},
\end{equation}
with $\tau'=\tau'(p,\tau):=\beta-\theta=(N'-2) (p-1)-(4+\tau)>-2$ (due to
$p>\tilde p_c(N',\tau)> \frac{N'+\tau}{N'-2}$).

We now show that Theorem \ref{t4.3} can be used to conclude the proof. To this end, we
need to analyze the function $f(p)$ when $\tau$ is replaced by $\tau'$.
To stress the dependence of $f(p)$ on $\tau$, we
write
$f(p)=f_\tau(p)$. For $(p,\tau)$ given above, and $\tau'=\tau'(p,\tau)>-2$, we now consider
the function $f_{\tau'}(q)$ for $q\in (1,\infty)$. From our analysis on $f_\tau(p)$ we know that
\[
\begin{aligned}
&f_{\tau'}(q)<\frac{(N'-2)^2}{4} \;\;\forall q\in (1, \tilde p_c(N',\tau'))\cup (p_c(N',\tau'),\infty),\\
&f_{\tau'}(q)>\frac{(N'-2)^2}{4} \;\; \forall q\in (\tilde p_c(N',\tau'), p_c(N',\tau')).
\end{aligned}
\]
A simple calculation shows that
\[
p-\frac{N'+2+2\tau}{N'-2}=\frac{N'+2+2\tau'}{N'-2}-p
\; \mbox{
and
$f_{\tau'}(p)=f_\tau(p)$.}
\]
Thus under our assumption on $p$, we have
\[
p>\frac{N'+2+2\tau'}{N'-2} \mbox{ and }
f_{\tau'}(p)=f_\tau(p)>\frac{(N'-2)^2}{4}.
\]
By the property of the function $f_{\tau'}(q)$, the above inequalities imply that
\[
p\in \left(\frac{N'+2+2\tau'}{N'-2}, p_c(N',\tau')\right).
\]
 In view of \eqref{4.6}, we
conclude that
\[
\frac{N'+2+2\tau'}{N'-2}<p<\min\{p_c(N',\tau'), p_c(N,0)\}.
\]
Therefore Theorem \ref{t4.3} applies to \eqref{4.7}, and $w(y)$ has fast decay at
$\infty$. This  implies that $x=0$ is a removable singularity of
$v$.

Finally we show that $V_\infty$ is stable in $\R^N\setminus\{0\}$
when $\frac{N'+\tau}{N'-2}<p\leq \tilde p_c(N',\tau)$. This is
equivalent to showing that \eqref{add} holds for such $p$, which
would follow if
\[
\frac{(N'-2)^2}{4}-p C_0^{p-1}=\frac{(N'-2)^2}{4}-f(p) \geq 0.
\]
But for $p\in (\frac{N'+\tau}{N'-2},\tilde p_c(N',\tau)]$ we do have
$f(p)\leq \frac{(N'-2)^2}{4}$. Thus $V_\infty$ is indeed stable
for such $p$. \qed

\end{document}